\documentclass{amsart}
\usepackage{amscd,amssymb}
\usepackage[all]{xy}
\xyoption{poly}

\newcommand{\gr}{{\operatorname{gr}\nolimits}}

\renewcommand{\Im}{\operatorname{Im}\nolimits}
\newcommand{\Ker}{\operatorname{Ker}\nolimits}

\newcommand{\rad}{\operatorname{rad}\nolimits}
\newcommand{\rrad}{\mathfrak{r}}

\newcommand{\Ext}{\operatorname{Ext}\nolimits}

\newcommand{\op}{{\operatorname{op}\nolimits}}

\newcommand{\G}{\Gamma}
\renewcommand{\L}{\Lambda}

\newcommand{\F}{{\mathcal F}}

\newcommand{\E}{{\mathcal E}}
\newcommand{\I}{{\mathcal I}}
\newcommand{\R}{{\mathcal R}}

\newcommand{\HH}{\operatorname{HH}\nolimits}

\newcommand{\rl}{{\operatorname{rl}\nolimits}}
\newcommand{\mo}{\mathfrak{o}}
\newcommand{\mt}{\mathfrak{t}}

\renewcommand{\G}{{\mathcal G}}
\newcommand{\J}{{\mathcal J}}
\newcommand{\N}{{\mathcal N}}

\newtheorem{lem}{Lemma}[section]
\newtheorem{prop}[lem]{Proposition}
\newtheorem{cor}[lem]{Corollary}
\newtheorem{thm}[lem]{Theorem}
\theoremstyle{definition}
\newtheorem{defin}[lem]{Definition}
\newtheorem*{remark}{Remark}
\newtheorem{example}[lem]{Example}
\begin{document}

\title[The Hochschild cohomology ring \dots]
{The Hochschild cohomology ring modulo nilpotence of a monomial algebra}
\author[Green]{Edward L. Green}
\address{Edward L. Green\\
Department of Mathematics\\ Virginia Tech\\
Blacksburg, VA 24061--0123\\
USA}
\email{green@math.vt.edu}
\author[Snashall]{Nicole Snashall}
\address{Nicole Snashall\\ Department of Mathematics\\
University of Leicester\\
University Road\\
Leicester, LE1 7RH\\
England}
\email{N.Snashall@mcs.le.ac.uk}
\author[Solberg]{\O yvind Solberg}
\address{\O yvind Solberg\\Institutt for matematiske fag\\
NTNU\\ N--7491 Trondheim\\ Norway}
\email{oyvinso@math.ntnu.no}

\begin{abstract}
For a finite dimensional monomial algebra $\L$ over a field $K$ we
show that the Hochschild cohomology ring of $\L$ modulo the ideal
generated by homogeneous nilpotent elements is a commutative finitely
generated $K$-algebra of Krull dimension at most one. This was
conjectured to be true for any finite dimensional algebra over a field
in \cite{SS}.
\end{abstract}

\dedicatory{Dedicated to Claus M. Ringel on the occasion of his
sixtieth birthday}

\date{\today}
\maketitle

\section*{Introduction}

Let $\L$ be a finite dimensional algebra over a field $K$, with
Jacobson radical $\rrad$. Denote by $\L^e$ the enveloping algebra
$\L^\op\otimes_K \L$ of $\L$. The Hochschild cohomology ring
$\HH^*(\L)$ of $\L$ is given by $\HH^*(\L)=\oplus_{i\geq
0}\Ext^i_{\L^e}(\L,\L)$ with the Yoneda product. It is well-known that
$\HH^*(\L)$ is a graded commutative ring, so that when the
characteristic of $K$ is different from two, then every homogeneous
element of odd degree squares to zero. Hence $\HH^*(\L)/\N$ is a
commutative $K$-algebra, where $\N$ is the ideal generated by the
homogeneous nilpotent elements. In \cite{SS} it was conjectured that
$\HH^*(\L)/\N$ is a finitely generated $K$-algebra. This paper is
devoted to showing this for any finite dimensional monomial
$K$-algebra $\L$ and in addition that $\HH^*(\L)/\N$ has Krull
dimension at most one. Recall that $\L$ is such a monomial algebra if
$\L=K\mathcal{Q}/I$ for some quiver $\mathcal{Q}$ and an ideal $I$ in
$K\mathcal{Q}$ generated by monomials of length at least two.

This conjecture was earlier known to be true for Nakayama algebras
with one relation (\cite{SS}), finite dimensional selfinjective
indecomposable algebras of finite representation type over an
algebraically closed field (\cite{GSS}), any block of a group ring of
a finite group (\cite{E,Go,V}), and any block of a finite dimensional
cocommutative Hopf algebra (\cite{FS} and see \cite{SS}). As a special
case of our results we point out the following consequence of the work
in \cite{C}. Namely, the conjecture is true for a path algebra
$K\mathcal{Q}/J^2$, where $J$ is the ideal generated by all the
arrows.

Even though the proofs are technical, the strategy of the proof is
easy to explain. The functor $\L/\rrad\otimes_\L -$ from right
$\L^e$-modules to right $\L$-modules induces a homomorphism of graded
rings $\varphi_{\L/\rrad}\colon \HH^*(\L)\to E(\L)$, where
$E(\L)=\oplus_{i\geq 0}\Ext^i_\L(\L/\rrad,\L/\rrad)$. It is shown in
\cite{SS} that $\Im \varphi_{\L/\rrad}$ is contained in the graded
centre $Z_\gr(E(\L))$ of $E(\L)$, that is, the subring of $E(\L)$
generated by all homogeneous elements $z$ such that
$zg=(-1)^{|z||g|}gz$ for each homogeneous element $g$ in $E(\L)$ where
$|x|$ denotes the degree of a homogeneous element $x$.  Furthermore the
map $\varphi_{\L/\rrad}$ induces an inclusion
$\overline{\varphi_{\L/\rrad}}$ of $\HH^*(\L)/\N$ into
$Z_\gr(E(\L))/\N_Z$, where $\N_Z$ is the ideal in $Z_\gr(E(\L))$
generated by homogeneous nilpotent elements. The method of proof is
then as follows. First we show that $Z_\gr(E(\L))/\N_Z$ is a
commutative finitely generated $K$-algebra, and secondly that beyond
some degree the even parts of the graded rings $\HH^*(\L)/\N$ and
$Z_\gr(E(\L))/\N_Z$ are isomorphic via the map
$\overline{\varphi_{\L/\rrad}}$. It then follows easily that
$\HH^*(\L)/\N$ is a finitely generated $K$-algebra.

The precise relationship between $\Im \varphi_{\L/\rrad}$ and $E(\L)$
for any finite dimensional $K$-algebra $\L$ has been found by Bernhard
Keller \cite[Theorem 3.5]{K}. The image $\Im \varphi_{\L/\rrad}$ is
shown to be the $\mathbb{A}_\infty$-centre of $E(\L)$. For further
details see section 7.

The paper is organized as follows. Let $\L=K\mathcal{Q}/I$ be a
monomial algebra. The first section is devoted to giving the
multiplicative basis of $E(\L)$ and the properties of this basis 
which are used throughout this paper. This is heavily based on the
papers \cite{Ba}, \cite{GHZ} and \cite{GZ}. The basis elements in $E(\L)$
correspond to paths in $K\mathcal{Q}$. In section three it is shown
that elements in $Z_\gr(E(\L))$ are obtained from closed walks and
equivalence classes of such. Non-nilpotent elements in $Z_\gr(E(\L))$
of even degree are made up of closed walks which are a product of
relations, and we call them tightly packed, extending walks.
Definitions and properties of these closed walks are given in section
two. The structure of non-nilpotent homogeneous elements in
$Z_\gr(E(\L))$ is analysed in further detail in the fourth section.
The next section is devoted to characterizing the elements $z^m$ for
homogeneous non-nilpotent elements $z$ in $Z_\gr(E(\L))$ for an even
positive integer $m$ less than or equal to the radical length of $\L$. In
the two last sections $Z_\gr(E(\L))$ and $\HH^*(\L)$ modulo the ideal
generated by homogeneous nilpotent elements are shown to be finitely
generated (commutative) $K$-algebras of Krull dimension at most one.
The paper ends with a brief discussion of the relationship between
$\Im \varphi_{\L/\rrad}$ and $E(\L)$ for a general finite dimensional
$K$-algebra and some examples illustrating our results.

\section*{Acknowledgement}
The authors acknowledge support from NSA and the Research Council of
Norway, and support and hospitality from the departments of
mathematics at the University of Leicester, Norwegian University of Science
and Technology and Virginia Tech.

\section{Preliminaries}

In \cite{GZ}, it was shown that there is a multiplicative basis for
the Ext algebra $E(\L)$ for a finite dimensional monomial algebra $\L
= K{\mathcal Q}/I$. We now describe this basis in terms of the paths
in $K{\mathcal Q}$. Arrows in a path are read from left to right. An
arrow $a$ starts at the vertex $\mathfrak{o}(a)$ and ends at the
vertex $\mathfrak{t}(a)$. We say that a path $p$ starts at the vertex
$\mathfrak{o}(p)$ and ends at the vertex $\mathfrak{t}(p)$. We denote
the {\it length} of a path $p$ by $\ell(p)$. A path $p$ is a {\it
  prefix} of a path $q$ if there is some path $p'$ with $q = pp'$. A
path $p$ is a {\it suffix} of a path $q$ if there is some path $p'$
with $q = p'p$.

Fix some minimal generating set of monomials for the ideal $I$ and
denote this set by $\R^2$.

We start by recalling some definitions and results from \cite{Ba},
\cite{GHZ} and \cite{GZ} about overlaps of paths in $K{\mathcal Q}$.

\bigskip

An overlap relation may be illustrated in the following way, where $P,
Q, U, V$ are paths in $K{\mathcal Q}$ with $VQ = PU$. 
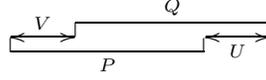
\[\xymatrix@W=0pt@M=0.3pt{%
\ar@{^{|}-^{|}}@<-1.25ex>[rrr]_P\ar@{{<}-{>}}[r]^{V} &
\ar@{_{|}-_{|}}@<1.25ex>[rrr]^{Q} & & \ar@{{<}-{>}}[r]_{U} & }\] %

We now give the precise definition of an overlap relation.

\begin{defin}[\cite{GHZ}]\label{defin:overlaprelation}
An {\it overlap relation} is a tuple $(P, Q, U, V)$ of paths in 
$K{\mathcal Q}$ such that
\begin{enumerate}
\item[(i)] $PU = VQ$, and
\item[(ii)] $1 \leq \ell(U) < \ell(Q)$ and $1 \leq \ell(V) < \ell(P)$.
\end{enumerate}
In this case we say $P$ {\it is overlapped by} $Q$, or that $Q$
{\it overlaps} $P$.
\end{defin}

Before giving the formal definition of a left overlap sequence and a
right overlap sequence we again explain the concepts pictorially and
by example.

A left overlap sequence is a sequence $(R_2, R_3, \ldots , R_n)$ with $R_2,
R_3, \ldots , R_n \in \R^2$ that, in the case $n=5$, may be illustrated
thus.
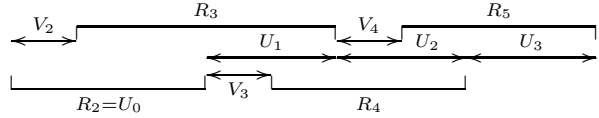
\[\xymatrix@W=0pt@M=0.3pt{%
\ar@{^{|}-^{|}}@<-2.75ex>[rrr]_{R_2=U_0}\ar@{<->}@<1.5ex>[r]^{V_2} &
\ar@{_{|}-_{|}}@<2.75ex>[rrrr]^{R_3} & & 
\ar@{<->}@<-1.5ex>[r]_{V_3}\ar@{<->}[rr]^{U_1} &
\ar@{^{|}-^{|}}@<-2.75ex>[rrr]_{R_4} &
\ar@{<->}[rr]^(0.7){U_2}\ar@{<->}@<1.5ex>[r]^{V_4} &
\ar@{_{|}-_{|}}@<2.75ex>[rrr]^{R_5} & 
\ar@{<->}[rr]^{U_3} & & }\]
The underlying path is $U_0U_1U_2U_3$ and we have $U_0U_1 = V_2R_3, U_1U_2 =
V_3R_4$ and $U_2U_3 = V_4R_5$.

A right overlap sequence is a sequence $(R_2, R_3, \ldots , R_n)$
with $R_2, R_3, \ldots , R_n \in \R^2$ that, in the case $n=5$, 
may be shown in the following way.
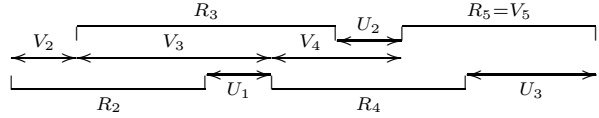
\[\xymatrix@W=0pt@M=0.3pt{%
\ar@{^{|}-^{|}}@<-2.75ex>[rrr]_{R_2}\ar@{<->}[r]^{V_2} &
\ar@{<->}[rrr]^{V_3}\ar@{_{|}-_{|}}@<2.75ex>[rrrr]^{R_3} & &
\ar@{<->}@<-1.5ex>[r]_{U_1} & 
\ar@{^{|}-^{|}}@<-2.75ex>[rrr]_{R_4}\ar@{<->}[rr]^(0.3){V_4} &
\ar@{<->}@<1.5ex>[r]^{U_2} & 
\ar@{_{|}-_{|}}@<2.75ex>[rrr]^{R_5=V_5} & \ar@{<->}@<-1.5ex>[rr]_{U_3}
& & }\]
The underlying path is $V_2V_3V_4V_5$ and we have $V_2V_3 = R_2U_1, V_3V_4 =
R_3U_2$ and $V_4V_5 = R_4U_3$.

\begin{example}\label{ex:overlaps}
\newcommand{\udspace}{\rule{0mm}{2.5mm}\raisebox{-0.8mm}{\rule{0mm}{0.8mm}}}
Let $Q$ be the quiver given by 
$$\xymatrix{1\ar[r]^a & 2\ar[r]^b & 3\ar[r]^c & 
4\ar[r]^d & 5\ar[r]^e & 6\ar[r]^f & 7\ar[r]^g & 8 }$$ 
and relations $\R^2 = \{ abc,bcde,cdef,efg\}$. 
Then $(abc,cdef)$ is a left overlap sequence with path
$\xymatrix@W=0pt@M=0pt@C=0pt{%
\ar@{^{|}-^{|}}@<-1.8ex>[rrrr] & \udspace a & \udspace b
\ar@{_{|}-_{|}}@<1.8ex>[rrrrr] & \udspace c & \udspace d &
\udspace e & \udspace f & }$. The pair $(abc,bcde)$ is also a left overlap
sequence and has path
$\xymatrix@W=0pt@M=0pt@C=0pt{%
\ar@{^{|}-^{|}}@<-1.8ex>[rrrr] & \udspace a
\ar@{_{|}-_{|}}@<1.8ex>[rrrrr] & \udspace b & \udspace c &
\udspace d & \udspace e & }$.
\end{example}

We now give the formal definitions of these concepts.

\begin{defin}[\cite{GHZ}] Let $n \geq 2$ and let $R_2 , R_3, \ldots , R_n \in
\R^2$.
\begin{enumerate}
\item The sequence $(R_2, R_3, \ldots , R_n)$ is a {\it left
overlap sequence} with data $$(U_0, \ldots , U_{n-2}, V_2, \ldots ,
V_{n-1})$$ if
 \begin{enumerate}
 \item[(i)] $U_0 = R_2$, and \item[(ii)] $(U_i, R_{i+3}, U_{i+1},
 V_{i+2})$ is an overlap relation for $0 \leq i \leq n-3$ (and in
 particular $U_iU_{i+1} = V_{i+2}R_{i+3}$).
 \end{enumerate}
The {\it path associated to a left overlap sequence} $(R_2, R_3, \ldots
, R_n)$ is the path $U_0\cdots U_{n-2}$.
\item The sequence $(R_2, R_3, \ldots , R_n)$ is a {\it right
overlap sequence} with data $$(U_1, \ldots , U_{n-2}, V_2, \ldots ,
V_n)$$ if
 \begin{enumerate}
 \item[(i)] $V_n = R_n$, and \item[(ii)] $(R_i, V_{i+1}, U_{i-1},
 V_i)$ is an overlap relation for $2 \leq i \leq n-1$ (and in
 particular $V_iV_{i+1} = R_iU_{i-1}$).
 \end{enumerate}
The {\it path associated to a right overlap sequence} $(R_2, R_3, 
\ldots , R_n)$ is the path $V_2\cdots V_n$.
\end{enumerate}
\end{defin}

\begin{remark}
(1) Note that the paths $P$ and $Q$ in Definition
\ref{defin:overlaprelation} are not required to be in $\R^2$. However all the
terms $R_2, R_3, \ldots , R_n$ in the definitions of a left overlap sequence and
a right overlap sequence must lie in $\R^2$.

(2) If $R_2$ and $R_3$ are in $\R^2$ so that $R_3$ overlaps $R_2$,
then $(R_2, R_3)$ is both a left overlap sequence and a right overlap
sequence.
\end{remark}

We now consider the algebra of Example \ref{ex:overlaps} 
again to discuss maximal left overlap sequences
and maximal right overlap sequences.

\begin{example}\label{ex2:overlaps}
\newcommand{\udspace}{\rule{0mm}{2.5mm}\raisebox{-0.8mm}{\rule{0mm}{0.8mm}}}
Let $Q$ be the quiver given by \[\xymatrix{1\ar[r]^a & 2\ar[r]^b &
3\ar[r]^c & 4\ar[r]^d & 5\ar[r]^e & 6\ar[r]^f & 7\ar[r]^g & 8 }\]
and relations $\R^2 = \{ abc,bcde,cdef,efg\}$. The left overlap sequence 
$(abc,cdef)$ with path
$\xymatrix@W=0pt@M=0pt@C=0pt{%
\ar@{^{|}-^{|}}@<-1.8ex>[rrrr] & \udspace a & \udspace b
\ar@{_{|}-_{|}}@<1.8ex>[rrrrr] & \udspace c & \udspace d &
\udspace e & \udspace f & }$
is not a maximal left overlap sequence. However
$(abc,bcde)$ with path
$\xymatrix@W=0pt@M=0pt@C=0pt{%
\ar@{^{|}-^{|}}@<-1.8ex>[rrrr] & \udspace a
\ar@{_{|}-_{|}}@<1.8ex>[rrrrr] & \udspace b & \udspace c &
\udspace d & \udspace e & }$ is a maximal left overlap sequence.
Also $(abc, bcde, efg)$, illustrated as
$\xymatrix@W=0pt@M=0pt@C=0pt{%
\ar@{^{|}-^{|}}@<-1.8ex>[rrrr] & \udspace a
\ar@{_{|}-_{|}}@<1.8ex>[rrrrr] & \udspace b & \udspace c &
\udspace d\ar@{^{|}-^{|}}@<-1.8ex>[rrrr] & \udspace e & \udspace f
& \udspace g & }$, is a maximal left overlap sequence with path
$abcdefg$. This path $abcdefg$ 
%is in $\mathcal{R}^4$, and 
is also the path associated to the maximal right overlap sequence
$(abc, cdef, efg)$ or 
$\xymatrix@W=0pt@M=0pt@C=0pt{%
\ar@{^{|}-^{|}}@<-1.8ex>[rrrr] & \udspace a
 & \ar@{_{|}-_{|}}@<1.8ex>[rrrrr]
\udspace b & \udspace c & \udspace d\ar@{^{|}-^{|}}@<-1.8ex>[rrrr]
& \udspace e & \udspace f & \udspace g & }$.
\end{example}

The formal definitions are as follows.

\begin{defin}[\cite{GHZ}]\label{defin:maximalos} 
Let $n \geq 2$ and let $R_2, R_3, \ldots , R_n \in
\R^2$. 
\begin{enumerate}
\item The left overlap sequence $(R_2, R_3, \ldots ,
R_n)$ is a {\it maximal left overlap sequence} if
  \begin{enumerate}
  \item[(i)] $V_2R_3 \neq PRQ$ for all nontrivial paths $P, Q$ and all $R \in
  \R^2$, and
  \item[(ii)] for $i = 3, \ldots , n-1$, $V_iR_{i+1} \neq PRQ$
  for all nontrivial paths $Q$ and all $R \in \R^2$.
  \end{enumerate}
\item The right overlap sequence $(R_2, R_3, \ldots ,
R_n)$ is a {\it maximal right overlap sequence} if
  \begin{enumerate}
  \item[(i)] $R_{n-1}U_{n-2} \neq PRQ$ for all nontrivial paths $P, Q$ and all $R \in
  \R^2$, and
  \item[(ii)] for $i = 1, \ldots , n-3$, $R_{i+1}U_i \neq PRQ$
  for all nontrivial paths $P$ and all $R \in \R^2$.
  \end{enumerate}
\end{enumerate}
In both cases we say the overlap sequence $(R_2, R_3, \ldots , R_n)$ is
of {\it length} $n$.
\end{defin}

For a monomial algebra, the paths associated to maximal left and maximal right
overlap sequences are closely related as we recall from \cite{Ba}.

\begin{thm}[\cite{Ba}]\label{thm:leftright}
For a monomial algebra, the path associated to a maximal left
overlap sequence $(R_2, R_3, \ldots , R_n)$ is the path associated to
some maximal right overlap sequence $(\tilde{R}_2, \tilde{R}_3, \ldots ,
\tilde{R}_n)$, and vice versa.
\end{thm}

\begin{remark}
It follows from the definitions that, in Theorem
\ref{thm:leftright}, we have $R_2 = \tilde{R}_2$ and $R_n =
\tilde{R}_n$. However $R_i$ need not equal $\tilde{R}_i$, for
$i=3, \ldots , n-1$, as Example \ref{ex2:overlaps} illustrates.
\end{remark}

The maximal left and right overlap sequences shall be omnipresent in
all of our arguments, so that we introduce appropriate notation for
these sets. 

\begin{defin}
Define $\R^n$ to be the set of paths associated to the maximal
left (or equivalently right) overlap sequences of length $n$ for $n\geq 2$.
This is consistent with the definition already introduced of the set $\R^2$. 
Let $\R^0$ denote the set of vertices of $\mathcal Q$ and let $\R^1$ denote 
the set of arrows of $\mathcal Q$.

From now on we use lower case letters, $r, r_i$, to denote the elements of
$\R^2$. By the term {\it relation}, we mean an element of the set $\R^2$,
that is, of our fixed minimal set of monomial generators for the ideal $I$. 
Let $(r_2, r_3, \ldots , r_n)$ be a maximal left (or equivalently right) overlap
sequence. We say that $r_2$ is the {\it first relation} and that $r_n$ is
the {\it last relation} of the path corresponding to this overlap sequence.
\end{defin}

The following are elementary properties of overlap sequences, and the
proofs are straightforward and left to the reader. 

\begin{lem}\label{lem:overlaps}
\begin{enumerate}
\item[(1)] The pair $(r_2,r_3)$ is a left overlap sequence if and
only if $(r_2,r_3)$ is a right overlap sequence.

\item[(2)] If $(r_2,r_3)$ is a left overlap sequence, then there
exists a maximal left overlap sequence $(r_2,\tilde{r}_3)$ such
that the path associated to $(r_2,\tilde{r}_3)$ is a prefix of the
path associated to $(r_2,r_3)$.

\item[(2$'$)] If $(r_2,r_3)$ is a right overlap sequence, then
there exists a maximal right overlap sequence $(\tilde{r}_2,r_3)$
such that the path associated to $(\tilde{r}_2,r_3)$ is a suffix
of the path associated to $(r_2,r_3)$.

\item[(3)] \sloppy Suppose $S=(r_2,r_3,\ldots,r_n)$ is a maximal left
overlap sequence and $S'=(r_2,r_3,\ldots,r_n,r_{n+1})$ is a left
overlap sequence. Then there exists
$\tilde{S}=(r_2,r_3,\ldots,r_n,\tilde{r}_{n+1})$ which is a
maximal left overlap sequence and the path associated to
$\tilde{S}$ is a prefix of the path associated to $S'$.

\item[(3$'$)] \sloppy Suppose $S=(r_3,r_4,\ldots,r_n)$ is a maximal
right overlap sequence and $S'=(r_2,r_3,r_4, \ldots,r_n)$ is a
right overlap sequence. Then there exists
$\tilde{S}=(\tilde{r}_2,r_3,r_4\ldots,r_n)$ which is a maximal
right overlap sequence and the path associated to $\widetilde{S}$
is a suffix of the path associated to $S'$.

\item[(4)] If $S=(r_2,r_3,\ldots,r_n)$ is a left overlap sequence,
then for all $i$ and $j$ with $1\leq i < j \leq n-1$, we have that
$S'=(r_{i+1},r_{i+2},\ldots,r_{j+1})$ is a left overlap sequence.
Note that even if $S$ is maximal we need not have $S'$ maximal.
However $S'$ is maximal in the case when $S$ is maximal and $i=1$.

\item[(4$'$)] If $S=(r_2,r_3,\ldots,r_n)$ is a right overlap
sequence, then for all $i$ and $j$ with $1\leq i < j \leq n-1$, we
have that $S'=(r_{i+1},r_{i+2},\ldots,r_{j+1})$ is a right overlap
sequence. Note that even if $S$ is maximal we need not have $S'$
maximal. However $S'$ is maximal in the case when $S$ is maximal
and $j=n-1$.
\end{enumerate}
\end{lem}

The importance of the sets $\R^n$ comes from the fact that they
describe a basis for the $n$-th graded part of the graded algebra
$E(\L)=\oplus_{i\geq 0}\Ext^i_\L(\L/\rrad,\L/\rrad)$. We explain this
next. Following \cite{GHZ}, fix a minimal $\L$-projective resolution
$(P^*, d^*)$ of the right $\L$-module $\L/\rrad$ which is determined
by $\bigcup_{m \geq 0}\R^m$. Then
$$P^m = \coprod_{R_j^m \in \R^m} \mt(R_j^m)\L .$$
For each $R_i^m \in \R^m$ there is a corresponding
element $g_i^m$ in $\Ext_{\L}^m(\L/\rrad, \L/\rrad)$. This element
$g_i^m$ is represented by the map
$$P^m \to \L/\rrad$$
given by
$$\mt(R_j^m) \mapsto \left\{
\begin{array}{ll}
\mt(R_i^m) & \mbox{if $j = i$,}\\
0 & \mbox{otherwise.}
\end{array}\right.$$

The set $\{g_i^m \colon R_i^m \in \R^m\}$ is a basis for
$\Ext_{\L}^m(\L/\rrad, \L/\rrad)$ and is denoted by $\G^m$. From
\cite{GZ}, the set $\bigcup_{m\geq 0}\G^m$ forms a multiplicative
basis for $\Ext_{\L}^*(\L/\rrad, \L/\rrad)$, that is, for $g_i^m
\in \G^m$ and $g_j^n \in \G^n$, if $g_i^mg_j^n \neq 0$ then
$g_i^mg_j^n \in \G^{m+n}$. For ease of notation and where no
confusion will arise, we may denote an element of $\G^m$ by $g^m$
rather than $g_i^m$ for some $i$; in this case the corresponding
path in $\R^m$ is denoted $R^m$. All relations between basis
elements are generated by
\begin{enumerate}
\item If $m = \sum_{i=1}^sm_i$ and if $R^{m_1}\cdots R^{m_s}$ is not in
$\R^m$, 
then $g^{m_1}\cdots g^{m_s} = 0$.
\item If $$R^{m_1}\cdots R^{m_s} = R^{n_1}\cdots R^{n_t}$$ is an element in
$\R^m$ where $m = \sum_{i=1}^s m_i = \sum_{j=1}^t n_j$, then 
$$g^{m_1}\cdots g^{m_s} = g^{n_1}\cdots g^{n_t}.$$
\end{enumerate}
Note that if $m = \sum_{i=1}^sm_i$ and if $R^{m_1}\cdots R^{m_s} \in \R^m$, 
then $g^{m_1}\cdots g^{m_s} \in \G^m$ and so
$g^{m_1}\cdots g^{m_s} \neq 0$.

In this paper we work extensively with the sets $\R^n$. Thus, for
$R^m_i$ in $\R^m$ and $R^n_j$ in $\R^n$ with corresponding elements
$g^m_i$ in $\G^m$ and $g^n_j$ in $\G^n$, we have $g^m_ig^n_j\neq 0$ if
and only if $R^m_iR^n_j$ is in $\R^{m+n}$. 
\bigskip

{F}rom the construction of this basis we have the following facts about the
elements $g_i^m$ and the corresponding paths $R_i^m$.  These results will be
used throughout this paper.

\begin{prop}[\cite{Ba,GHZ,GZ}]\label{prop:extcancel}
\sloppy For a finite dimensional monomial algebra $\L = K{\mathcal Q}/I$,
the following properties hold.
\begin{enumerate}
\item If $R_i^uR_j^v = R_i^uR_k^w \in \R^{u+v}$ then $v = w$ and
$R_j^v = R_k^w$.
\item If $R_j^vR_i^u = R_k^wR_i^u \in \R^{u+v}$ then $v = w$ and
$R_j^v = R_k^w$.
\item If $g_i^ug_j^v = g_i^ug_k^w \neq 0$ then $v = w$ and $g_j^v = g_k^w$.
\item If $g_j^vg_i^u = g_k^wg_i^u \neq 0$ then $v = w$ and $g_j^v = g_k^w$.
\end{enumerate}
\end{prop}

\begin{prop}[\cite{GZ}]\label{prop:equal}
For a finite dimensional monomial algebra $\L = K{\mathcal Q}/I$,
if $g_1^{m_1}\ldots g_s^{m_s} = g_{1'}^{m_1}\ldots g_{s'}^{m_s} \neq 0$
then $g_i^{m_i} = g_{i'}^{m_i}$ for $i = 1, \ldots , s$.
\end{prop}

The next two results give further properties of the basis for $E(\L)$,
which we use later.

\begin{prop}\label{prop:uniquepaths}
For a finite dimensional monomial algebra $\L = K{\mathcal Q}/I$,
the following properties hold.
\begin{enumerate}
\item If $R_i^mp = R_j^mq$ for some paths $p$ and $q$, then $R_i^m = R_j^m$
and $p = q$.
\item If $pR_i^m = qR_j^m$ for some paths $p$ and $q$, then $R_i^m = R_j^m$
and $p = q$.
\end{enumerate}
\end{prop}

\begin{proof}
(1) Suppose $R_i^mp = R_j^mq$ but that $R_i^m \neq R_j^m$. Let
$R_i^m$ have maximal left overlap sequence $(r_2, r_3, \ldots ,
r_m)$ and let $R_j^m$ have maximal left overlap sequence
$(\tilde{r}_2, \tilde{r}_3, \ldots , \tilde{r}_m)$. Then by
construction, $r_2 = \tilde{r}_2$, and there is some $k$
such that $r_2 = \tilde{r}_2, \ldots , r_k = \tilde{r}_k$ but
$r_{k+1} \neq \tilde{r}_{k+1}$. Without loss of generality we may
assume that $\tilde{r}_{k+1}$ starts after $r_{k+1}$. Thus we have
the following diagram.
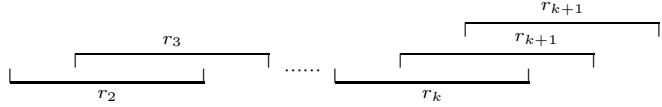
\[\xymatrix@W=0pt@M=0.3pt{%
\ar@{^{|}-^{|}}@<-1.25ex>[rrr]_{r_2} &
\ar@{_{|}-_{|}}@<1.25ex>[rrr]^{r_3} & &
\ar@{}@<-1.25ex>[rrr]^{\cdots\cdots} & &
\ar@{^{|}-^{|}}@<-1.25ex>[rrr]_{r_k} &
\ar@{_{|}-_{|}}@<1.25ex>[rrr]^(0.7){r_{k+1}} &
\ar@{_{|}-_{|}}@<4ex>[rrr]^{\tilde{r}_{k+1}} & & & }\] %
But this contradicts $(\tilde{r}_2, \tilde{r}_3, \ldots ,
\tilde{r}_m)$ being a maximal left overlap sequence. Hence $R_i^m
= R_j^m$ and thus $p=q$.

(2) The proof is similar.
\end{proof}

\begin{prop}\label{prop:multrelation}
For a finite dimensional monomial algebra $\L = K{\mathcal Q}/I$,
if $g^mg^n \neq 0$ and $n\geq 2$, then $g^mg^2 \neq 0$ where $g^2$
is the basis element corresponding to the first relation of $R^n$.
\end{prop}

\begin{proof}
Let $(r_2, r_3, \ldots , r_m, r_{m+1}, \ldots , r_{m+n})$ be a
maximal left overlap sequence for $R^mR^n$. Then from Lemma
\ref{lem:overlaps}, we have that $(r_2, r_3, \ldots , r_m)$ is a
maximal left overlap sequence for $R^m$. Since $\mt(R^m) =
\mo(R^n)$, the relation $r_{m+2}$ is the first relation of $R^n$.
Thus $(r_2, r_3, \ldots , r_m, r_{m+1}, r_{m+2})$ is a maximal
left overlap sequence for $R^mr_{m+2}$. Hence $R^mr_{m+2} \in
\R^{m+2}$ and so $g^mg^2 \neq 0$ where $g^2$ is the basis element
corresponding to the first relation of $R^n$.
\end{proof}

In multiplying the basis elements for $E(\L)$ the following concept is
central.

\begin{defin}[\cite{GHZ}]\label{def:tailpaths}
Let $\L = K{\mathcal Q}/I$ be a finite dimensional monomial algebra.
\begin{enumerate}
\item For each $R_i^m \in \R^m$, $m \geq 1$,
there exists a unique path $t$ and unique element $R_k^{m-1} \in
\R^{m-1}$ with $R_i^m = R_k^{m-1}t$. We call $t$ the {\it tail
path} of $R_i^m$ and write $t = t(R_i^m)$. 
\item For each $R_i^m \in \R^m$, $m \geq 1$, there exists a unique path
$b$ and unique element $R_j^{m-1} \in \R^{m-1}$ with $R_i^m =
bR_j^{m-1}$. We call $b$ the {\it beginning path} of $R_i^m$ and
write $b = b(R_i^m)$.
\end{enumerate}
\end{defin}

\bigskip

We end this section with some remarks about closed walks in a
quiver.

\begin{defin}
A {\it closed walk}
in $\mathcal{Q}$ is a non-trivial path $C$ in $K\mathcal{Q}$
such that $C = eCe$ for some vertex $e$.
\end{defin}

We do not make any assumptions with this terminology as to whether
or not $C$ is a nonzero element in the algebra $\L$.

\begin{lem}\label{lem:walk}
Let $p$, $q$ and $r$ be walks in a quiver ${\mathcal Q}$, and let
$a_1, \ldots , a_t$ be arrows in ${\mathcal Q}$. Suppose that
$pq=qr$. If $p=a_1a_2\cdots a_t$, then $q=p^sa_1\cdots a_i$ and
$r=a_{i+1}\cdots a_ta_1\cdots a_i$ for some $0 \leq i < t$ and $s\geq 0$,
where the case $i = 0$ implies that $q = p^s$ and $r = p$.
Moreover $p$ and $r$ are closed walks in ${\mathcal Q}$.
\end{lem}

\section{Tightly packed maximal overlap sequences}

Our strategy of proof is to first analyse the graded centre of $E(\L)$
modulo nilpotence. Recall that the graded centre $Z_\gr(E(\L))$ of
$E(\L)$ is the subring generated by all homogeneous elements $z$ in
$E(\L)$ such that $zg=(-1)^{|z||g|}gz$ for each homogeneous element
$g$ in $E(\L)$. The elements in $E(\L)$ are described by the basis
$\cup_{m\geq 0} \G^m$. The path $R^m$ corresponding to a basis vector
$g^m$ need not be a product of relations; however it turns out that all
the elements in $Z_\gr(E(\L))$ are associated to paths 
which are products of relations. Such paths we call
tightly packed. In this section we define and give some properties of
tightly packed paths which will be needed for our
characterization of the non-nilpotent elements of $Z_\gr(E(\L))$.

\begin{defin}\label{defin:tightlypacked}
A path $p$ in $K{\mathcal Q}$ is {\it tightly packed} if $p =
r_1r_2\cdots r_n$ for some $r_i \in \R^2$, $i=1, \ldots , n$.
We say that $r_1, r_2, \ldots , r_n$ are the {\it relations in} $p$.
\end{defin}

The next four results give basic properties of tightly packed paths.

\begin{prop}\label{prop:tplength}
Suppose $R^{2n}$ in $\R^{2n}$ is tightly packed with $R^{2n} =
r_2r_4\cdots r_{2m}$ for some $m$ and $r_2, r_4, \ldots , r_{2m}
\in \R^2$. Then $m=n$ and the maximal left (resp.\ right) overlap
sequence for $R^{2n}$ is of the form $(r_2, r_3, r_4, \ldots ,
r_{2n})$.
\end{prop}

\begin{proof}
Let $(\tilde{r}_2, \tilde{r}_3, \ldots , \tilde{r}_{2n})$ be the
maximal left overlap sequence for $R^{2n}$. Then the path
corresponding to this maximal left overlap sequence is of the form
$\tilde{r}_2p$ and hence $\tilde{r}_2 = r_2$. Now $\mt(r_2) =
\mo(r_4)$ so, by maximality of the left overlap sequence, we have
$\tilde{r}_4 = r_4$. By induction $\tilde{r}_6 = r_6, \ldots ,
\tilde{r}_{2n} = r_{2n}$. Then the path associated to the maximal
left overlap sequence is $r_2r_4\cdots r_{2n}$ and so $m=n$.
\end{proof}

\begin{prop}\label{prop:tpleftoverlap}
Suppose the path $p = r_2r_4\cdots r_{2n}$ is tightly packed and
that there is some left overlap sequence $(r_2, r_3, r_4, \ldots,
r_{2n})$. Then $p \in \R^{2n}$.
\end{prop}

\begin{proof}
\sloppy We show that there is some maximal left overlap sequence $(r_2,
\tilde{r}_3, r_4, \ldots , \tilde{r}_{2n-1}, r_{2n} )$ with path
$p$. Since $(r_2, r_3, r_4, \ldots, r_{2n})$ is a left overlap
sequence, the relation $r_3$ overlaps $r_2$. Thus there is some
maximal left overlap $\tilde{r}_3$ with $r_2$, and so
$\tilde{r}_3$ ends at or before $r_3$ ends. Now $\mt(r_2) = \mo(r_4)$ so
$(r_2, \tilde{r}_3, r_4)$ is a maximal left overlap sequence. The
relation $r_5$ overlaps $r_4$ and starts at or after the end of
$\tilde{r}_3$ as the diagram illustrates.
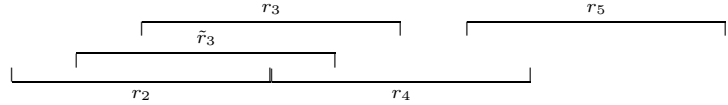
\[\xymatrix@W=0pt@M=0.3pt{%
\ar@{^{|}-^{|}}@<-1.25ex>[rrrr]_{r_2} &
\ar@{_{|}-_{|}}@<1.25ex>[rrrr]^{\tilde{r}_3} &
\ar@{_{|}-_{|}}@<4ex>[rrrr]^{r_3} & &
\ar@{^{|}-^{|}}@<-1.25ex>[rrrr]_{r_4} & & &
\ar@{_{|}-_{|}}@<4ex>[rrrr]^{r_5} & & & & }\] %
Thus there is some maximal overlap $\tilde{r}_5$ with $r_4$ such that
$\tilde{r}_5$ ends at or before $r_5$ ends and $(r_2, \tilde{r}_3, r_4,
\tilde{r}_5)$ is a maximal left overlap sequence. Continuing in this way we
have a maximal left overlap sequence $(r_2, \tilde{r}_3, r_4,
\tilde{r}_5, \ldots , \tilde{r}_{2n-1}, r_{2n})$ with path $p$.
Hence $p$ is in $\R^{2n}$.
\end{proof}

\begin{cor}\label{cor:tpsubpaths}
Suppose $R^{2n} = r_2r_4\cdots r_{2n}$ is tightly packed and
$R^{2n} \in \R^{2n}$. Then for each $1 \leq i \leq j \leq n$, the
subpath $r_{2i}r_{2i+2}\cdots r_{2j}$ is in $\R^{2(j-i)+2}$.
\end{cor}

\begin{proof}
\sloppy From Proposition \ref{prop:tpleftoverlap}, there is a maximal left
overlap sequence $(r_2, r_3, r_4, \ldots r_{2n})$ for $R^{2n}$.
From Lemma \ref{lem:overlaps}, $(r_{2i}, r_{2i+1}, \ldots ,
r_{2j})$ is a left overlap sequence for $r_{2i}r_{2i+2}\cdots
r_{2j}$. The result now follows by applying Proposition
\ref{prop:tpleftoverlap} again.
\end{proof}

\begin{prop}\label{prop:tppowers}
Suppose $R^{2n}$ is in $\R^{2n}$ and $(R^{2n})^u$ is in $\R^{2nu}$ for
some $u \geq 1$. Suppose also that $(R^{2n})^u$ is tightly packed.
Then $R^{2n}$ is tightly packed.
\end{prop}

\begin{proof}
Suppose $(R^{2n})^u$ in $\R^{2nu}$ is tightly packed. From Proposition
\ref{prop:tplength}, we may write $(R^{2n})^u = r_2r_4\cdots r_{2nu}$
with $r_2, r_4, \ldots r_{2nu} \in \R^2$. Then Corollary
\ref{cor:tpsubpaths} gives that $r_2r_4\cdots r_{2n}$ is in
$\R^{2n}$. Now, as paths, $(r_2r_4\cdots
r_{2n})(r_{2n+2}r_{2n+4}\cdots r_{2nu}) = (R^{2n})^u =
R^{2n}((R^{2n})^{(u-1)})$. So from Proposition
\ref{prop:uniquepaths} we have $r_2r_4\cdots r_{2n} = R^{2n}$.
Hence $R^{2n}$ is tightly packed.
\end{proof}

\begin{remark}
If $(R^{2n})^u \in \R^{2nu}$ for some $u\geq 2$, then $(R^{2n})^u$ is the
path corresponding to some element of $\R^{2nu}$, and in particular, $R^{2n}$
is a closed walk.
\end{remark}

When multiplying elements in $E(\L)$, the tail and the beginning of a
path corresponding to a basis element are important for deciding the
product of two basis elements. In particular we are interested in
non-nilpotent elements so that the tails of powers of basis elements
need to be considered. The next three results discuss the tails of
powers of tightly packed paths, culminating by
showing that the tails of high enough powers of tightly packed paths
become equal if the path is ``non-nilpotent''.

\begin{prop}\label{prop:tptailpath}
Suppose $R^{2n}$ in $\R^{2n}$ is tightly packed and $R^{2n}R^{2n}$
is in $\R^{4n}$. Then $t(R^{2n}R^{2n})$ is a suffix of $t(R^{2n})$.
\end{prop}

\begin{proof}
\sloppy Let $R^{2n}$ have maximal left overlap sequence
$(r_2,r_3,\ldots,r_{2n})$, and let
$(\tilde{r}_2, \tilde{r}_3,\ldots,\tilde{r}_{4n})$ be the maximal left overlap
sequence for $R^{2n}R^{2n}$.  Since $R^{2n}$ is tightly packed,
$R^{2n}=r_2r_4\cdots r_{2n}$. It follows that $R^{2n}R^{2n}$ is
tightly packed and so
$$r_2r_4\cdots r_{2n}r_2r_4\cdots r_{2n} =
R^{2n}R^{2n} = \tilde{r}_2\tilde{r}_4\cdots \tilde{r}_{4n}.$$ Thus
$\tilde{r}_i = \tilde{r}_{i+2n} = r_i$ for $i=2,4,6,\ldots,2n$.
Moreover, $(\tilde{r}_2,\ldots,\tilde{r}_{2n})$ is a maximal left
overlap sequence for $R^{2n}$ and so $\tilde{r}_2 = r_2,
\tilde{r}_3 = r_3, \ldots, \tilde{r}_{2n} = r_{2n}$.

We now consider $\tilde{r}_{2n+1}$. Viewing paths from left to
right we have the following maximal left overlap sequence
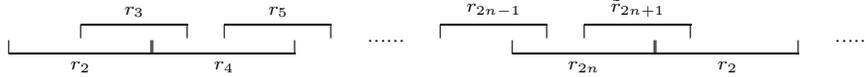
\[\xymatrix@W=0pt@M=0.3pt@C=13pt{%
\ar@{^{|}-^{|}}@<-1.25ex>[rrrr]_{r_2} & &
\ar@{_{|}-_{|}}@<1.25ex>[rrr]^{r_3} & &
\ar@{^{|}-^{|}}@<-1.25ex>[rrrr]_{r_4} & &
\ar@{_{|}-_{|}}@<1.25ex>[rrr]^{r_5} & & &  
\ar@{}@<-1.25ex>[rrr]^{\cdots\cdots} & & &
\ar@{_{|}-_{|}}@<1.25ex>[rrr]^{r_{2n-1}} & & 
\ar@{^{|}-^{|}}@<-1.25ex>[rrrr]_{r_{2n}} & & 
\ar@{_{|}-_{|}}@<1.25ex>[rrr]^{\tilde{r}_{2n+1}}& &
\ar@{^{|}-^{|}}@<-1.25ex>[rrrr]_{r_2} & & & & &  
\ar@{}@<-1.25ex>[r]^{\cdots\cdots} & }\] %
where either $\tilde{r}_{2n+1}$ does not end after $r_3$ begins or
$\tilde{r}_{2n+1}$ ends strictly after $r_3$ begins. In the first
case, if $\tilde{r}_{2n+1}$ does not end after $r_3$ begins, then
$\tilde{r}_{2n+1}$ is situated as follows:
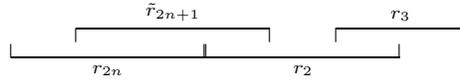
\[\xymatrix@W=0pt@M=0.3pt{%
\ar@{^{|}-^{|}}@<-1.25ex>[rrr]_{r_{2n}} & %
\ar@{_{|}-_{|}}@<1.25ex>[rrr]^{\tilde{r}_{2n+1}} & &
\ar@{^{|}-^{|}}@<-1.25ex>[rrr]_{r_2} & &
\ar@{_{|}-_{|}}@<1.25ex>[rr]^{r_3}& & }\] %
By maximality of the left overlap sequence
$(r_2,r_3,\ldots,r_{2n})$ for $R^{2n}$, the relation $r_3$ is the
first relation that overlaps $r_2$. Hence $\tilde{r}_{2n+3}=r_3$.
Then the maximal left overlap sequence for $R^{2n}R^{2n}$
continues as the sequence for $R^{2n}$, that is,
$\tilde{r}_{2n+i}=\tilde{r}_i$ for $i=1,2,\ldots,2n$. Thus
$t(R^{2n}R^{2n})=t(R^{2n})$.

So now consider the second case, that is, suppose
$\tilde{r}_{2n+1}$ ends strictly after $r_3$ begins so that
$\tilde{r}_{2n+1}$ is now situated as follows
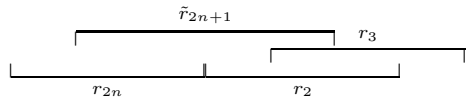
\[\xymatrix@W=0pt@M=0.3pt{%
\ar@{^{|}-^{|}}@<-1.25ex>[rrr]_{r_{2n}} &
\ar@{_{|}-_{|}}@<2.8ex>[rrrr]^{\tilde{r}_{2n+1}} & &
\ar@{^{|}-^{|}}@<-1.25ex>[rrr]_{r_2} &
\ar@{_{|}-_{|}}@<1.25ex>[rrr]^{r_3} & & & }\] %
Now $R^{2n}R^{2n}$ is in $\R^{4n}$ so there exists some relation
$\tilde{r}_{2n+3}$ along the path and this relation must overlap
$r_2$, ending strictly after $r_3$ ends, and thus we have
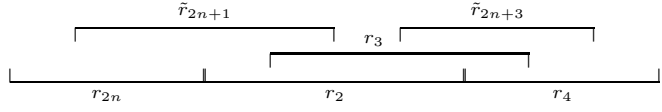
\[\xymatrix@W=0pt@M=0.3pt{%
\ar@{^{|}-^{|}}@<-1.25ex>[rrr]_{r_{2n}} &
\ar@{_{|}-_{|}}@<3.5ex>[rrrr]^{\tilde{r}_{2n+1}} & &
\ar@{^{|}-^{|}}@<-1.25ex>[rrrr]_{r_2} &
\ar@{_{|}-_{|}}@<1.25ex>[rrrr]^(.4){r_3} & &
\ar@{_{|}-_{|}}@<3.5ex>[rrr]^{\tilde{r}_{2n+3}} &
\ar@{^{|}-^{|}}@<-1.25ex>[rrr]_{r_4} & & & }\] %
Now either $\tilde{r}_{2n+3}$ does not end after $r_5$ begins, in which case
$t(R^{2n}R^{2n})=t(R^{2n})$ as above, or
$\tilde{r}_{2n+3}$ ends strictly after $r_5$ begins, in which case
we consider $\tilde{r}_{2n+5}$.

Continue inductively until we reach $\tilde{r}_{4n-3}$. Either
$\tilde{r}_{4n-3}$ does not end after $r_{2n-1}$ begins, in which
case $t(R^{2n}R^{2n})=t(R^{2n})$ as above, or $\tilde{r}_{4n-3}$
ends strictly after $r_{2n-1}$ begins, in which case
$\tilde{r}_{4n-3}$ is situated as follows:
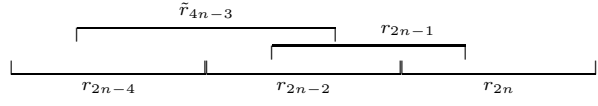
\[\xymatrix@W=0pt@M=0.3pt{%
\ar@{^{|}-^{|}}@<-1.25ex>[rrr]_{r_{2n-4}} &
\ar@{_{|}-_{|}}@<2.8ex>[rrrr]^{\tilde{r}_{4n-3}} & &
\ar@{^{|}-^{|}}@<-1.25ex>[rrr]_{r_{2n-2}} &
\ar@{_{|}-_{|}}@<1.25ex>[rrr]^(.7){r_{2n-1}} & &
\ar@{^{|}-^{|}}@<-1.25ex>[rrr]_{r_{2n}} & & & }\] %
Then $\tilde{r}_{4n-1}$ must end strictly after $r_{2n-1}$ ends.
Thus we have
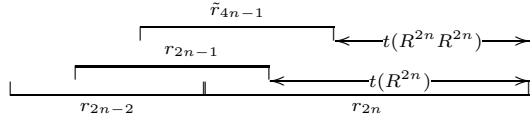
\[\xymatrix@W=0pt@M=0.3pt{%
\ar@{^{|}-^{|}}@<-1.25ex>[rrr]_{r_{2n-2}} &
\ar@{_{|}-_{|}}@<1.25ex>[rrr]^(0.6){r_{2n-1}} &
\ar@{_{|}-_{|}}@<4.7ex>[rrr]^{\tilde{r}_{4n-1}} &
\ar@{^{|}-^{|}}@<-1.25ex>[rrrrr]_{r_{2n}} &
\ar@{{<}-{>}}[rrrr]|{t(R^{2n})} & \ar@{{<}-{>}}@<3.5ex>[rrr]|{t(R^{2n}R^{2n})}%
 & & &
}\] %
and so $t(R^{2n}R^{2n})$ is a proper suffix of $t(R^{2n})$.
\end{proof}

\begin{remark}
The corresponding proof for the beginning paths shows that
$b(R^{2n}R^{2n})$ is a prefix of $b(R^{2n})$.
\end{remark}

\begin{prop}\sloppy 
Suppose $R^{2n}$ in $\R^{2n}$ is tightly packed and
$(R^{2n})^{(k+1)}$ is in $\R^{2n(k+1)}$ for some $k \geq 1$. Then
$t((R^{2n})^{(k+1)})$ is a suffix of $t((R^{2n})^k)$.
\end{prop}
\begin{proof}
Note that if $(R^{2n})^{(k+1)}$ is in $\R^{2n(k+1)}$, for some $k \geq 1$,
then $(R^{2n})^k$ is in $\R^{2nk}$ by Corollary \ref{cor:tpsubpaths}. 
The proof is now by induction with the case $k=1$ being Proposition
\ref{prop:tptailpath}. The rest of the proof is similar and is left to the
reader.
\end{proof}

\begin{defin} A closed walk $R^{2n}$ is {\it extending}
if $(R^{2n})^s \in \R^{2ns}$ for all $s\ge 1$.
\end{defin}

\begin{prop}\label{prop:tpalltails}
Suppose $R^{2n}$ is an extending closed walk. Then
there exists $u \geq 1$ such that
$t((R^{2n})^u)=t((R^{2n})^{(u+v)})$ for all $v\geq 1$. Moreover $u \leq
\rl(\L)-1$ where $\rl(\L)$ is the radical length of $\L$.
\end{prop}
\begin{proof}
From Proposition \ref{prop:tptailpath} the lengths of the tail paths
$t((R^{2n})^k)$ for $k\geq 1$ form a decreasing sequence, that is,
\[ \ell(t(R^{2n})) \geq \ell(t((R^{2n})^2)) \geq \ell(t((R^{2n})^3))
\geq \cdots .\] 
Since each successive tail path is a suffix of the previous tail
path, there exists some $u$ with $t((R^{2n})^u)=t((R^{2n})^{(u+1)})$.
Moreover each tail path has length at least 1, so $u \leq \rl(\L)-1$.

Let $(r_2,r_3,\ldots,r_{2nu})$ be a maximal left overlap sequence
for $(R^{2n})^u$, so that we have the following diagram:
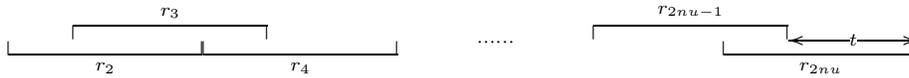
\[\xymatrix@W=0pt@M=0.3pt{%
\ar@{^{|}-^{|}}@<-1.25ex>[rrr]_{r_2} &
\ar@{_{|}-_{|}}@<1.25ex>[rrr]^{r_3} & &
\ar@{^{|}-^{|}}@<-1.25ex>[rrr]_{r_4} & & & &
\ar@{}@<-1.25ex>[r]^{\cdots\cdots} & &
\ar@{_{|}-_{|}}@<1.25ex>[rrr]^{r_{2nu-1}} & &
\ar@{^{|}-^{|}}@<-1.25ex>[rrr]_{r_{2nu}} &
\ar@{{<}-{>}}[rr]|{t} & & }\] %
Then there is a maximal left overlap sequence
$(r_2,r_3,\ldots,r_{2nu},\tilde{r}_{2nu+1},\ldots,\tilde{r}_{2n(u+1)})$
for $(R^{2n})^{(u+1)}$ with $\tilde{r}_{2nu+2}=r_2$,
$\tilde{r}_{2nu+4}=r_4$,\ldots, $\tilde{r}_{2n(u+1)}=r_{2n}$ since
$R^{2n}$ is tightly packed.

The relation $\tilde{r}_{2nu+1}$ overlaps $r_{2nu}$ and starts
after $r_{2nu-1}$ ends; we consider successively the overlaps
$\tilde{r}_{2nu+1}$, $\tilde{r}_{2nu+3}$, \ldots,
$\tilde{r}_{2n(u+1)-3}$, as in the proof of Proposition
\ref{prop:tptailpath}. Since $t((R^{2n})^{(u+1)})=t((R^{2n})^u)$, the relation
$\tilde{r}_{2n(u+1)-1}$ must end precisely where $r_{2nu-1}$ ends,
giving a maximal left overlap sequence
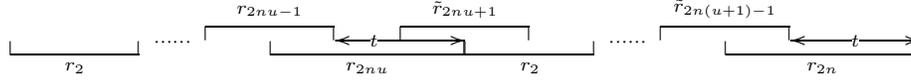
\[\xymatrix@W=0pt@M=0.3pt@R=6pt{%
\ar@{^{|}-^{|}}@<-1.25ex>[rr]_{r_2} & &
\ar@{}@<-1.25ex>[r]^{\cdots\cdots} &
\ar@{_{|}-_{|}}@<1.25ex>[rr]^{r_{2nu-1}} &
\ar@{^{|}-^{|}}@<-1.25ex>[rrr]_{r_{2nu}} &
\ar@{{<}-{>}}[rr]|(0.3){t} &
\ar@{_{|}-_{|}}@<1.25ex>[rr]^{\tilde{r}_{2nu+1}} &
\ar@{^{|}-^{|}}@<-1.25ex>[rr]_{r_2} & &
\ar@{}@<-1.25ex>[r]^{\cdots\cdots} &
\ar@{_{|}-_{|}}@<1.25ex>[rr]^{\tilde{r}_{2n(u+1)-1}} &
\ar@{^{|}-^{|}}@<-1.25ex>[rrr]_{r_{2n}} &
\ar@{{<}-{>}}[rr]|{t} & & }\] %
\sloppy Note that $\tilde{r}_{2n(u+1)-1}$ is necessarily equal to 
$r_{2nu-1}$. Then,
repeating the overlaps $\tilde{r}_{2nu+1}$, \ldots,
$\tilde{r}_{2n(u+1)-1}$, shows that
\[(r_2,r_3,\ldots,r_{2nu},\tilde{r}_{2nu+1},\ldots,
\tilde{r}_{2n(u+1)-1},\tilde{r}_{2nu+1},\ldots,
\tilde{r}_{2n(u+1)-1})\] is a maximal left overlap sequence for
$(R^{2n})^{(u+2)}$. Thus $t((R^{2n})^{(u+2)})=t((R^{2n})^{(u+1)})$. Hence, by
induction, $t((R^{2n})^u)=t((R^{2n})^{(u+v)})$ for all $v\geq 1$.
\end{proof}

As we have seen above, the tail paths of the powers of an extending closed
walk always become equal after the power $\rl(\L)-1$ of the walk. However we
could have equality for an earlier exponent. This gives rise to the next
definition, which is of the concept of being left or right stable.

\begin{defin} Let $R^{2n}\in \R^{2n}$ be an extending closed walk.
\begin{enumerate}
\item The path $R^{2n}$ is {\it left $j$-stable} if
  \begin{enumerate}
  \item $j$ is odd, and
  \item whenever $s$ is such that $2ns\geq j+1+2n$, then
   $r_j=r_{2n+j}$ and $r_{j+1}=r_{2n+j+1}$
  where $(r_2,r_3,\dots,r_{2ns})$ is the maximal left overlap sequence
  for $(R^{2n})^s$.
  \end{enumerate}
\item The path $R^{2n}$ is {\it right $j$-stable} if
  \begin{enumerate}
  \item $j$ is odd, and
  \item whenever $s$ is such that $2ns \geq j+1+2n$, then 
  $\tilde{r}_{2ns+2-j}=\tilde{r}_{2ns+2-j-2n}$ and 
  $\tilde{r}_{2ns+1-j}=\tilde{r}_{2ns+1-j-2n}$ 
  where
  $(\tilde{r}_2,\tilde{r}_3,\dots,\tilde{r}_{2ns})$ is the maximal right 
  overlap sequence for $(R^{2n})^s$.
  \end{enumerate}
\end{enumerate}
\end{defin}

The rest of this section is devoted to describing the tail path and the
beginning path of a power of a tightly packed, extending closed walk. 

\begin{prop}\label{conv:p1}
\begin{enumerate}
\item If $R^{2n}$ is a tightly packed, left $j$-stable, closed
walk with $j \leq 2n(s-1)-1$ and with maximal left overlap sequence
$(r_2,r_3,\dots,r_{2ns})$ for $(R^{2n})^s$,
then $r_k=r_{2n+k}$ for all $j \leq k \leq 2n(s-1)$.
\item If $R^{2n}$ is a tightly packed, right $j$-stable, closed
walk with maximal right overlap sequence
$(\tilde{r}_2,\tilde{r}_3,\dots,\tilde{r}_{2ns})$ for $(R^{2n})^s$, 
then $\tilde{r}_k=\tilde{r}_{k-2n}$ for all $2n+2 \leq k \leq j$.
\end{enumerate}
\end{prop}

\begin{proof} (1) \sloppy Since $R^{2n}$ is tightly packed, $r_{2i} = r_{2n+2i}$ for
$i = 1, \ldots , n(s-1)$. Let $j = 2i+1$ where $R^{2n}$ is left $j$-stable.
Then $r_{2n+2i+1} = r_{2i+1}, r_{2n+2i+2} = r_{2i+2}$ and $r_{2n+2i+4} =
r_{2i+4}$. By maximality of the left overlap sequence
$(r_2,r_3,\ldots,r_{2ns})$, it follows that
$r_{2i+3}=r_{2n+2i+3}$. Continuing inductively,
$r_{2i+5}=r_{2n+2i+5}, \ldots , r_{2n(s-1)-1}=r_{2ns-1}$. Hence
$r_k=r_{2n+k}$ for all $j \leq k \leq 2n(s-1)$.

(2) The proof is similar.
\end{proof}

\begin{prop}\label{conv:p2}
Let $R^{2n}\in\R^{2n}$ be a tightly packed, extending closed walk.  Then the
following statements hold.
\begin{enumerate}
\item The path $R^{2n}$ is left $j$-stable for some $j\ge 3$, and 
$t((R^{2n})^v) = t((R^{2n})^{(v+1)})$ for each $v$ such that $j<2nv$.
\item The path $R^{2n}$ is right $j$-stable for some $j\ge 3$, and 
$b((R^{2n})^v) = b((R^{2n})^{(v+1)})$ for each $v$ such that $j<2nv$.
\end{enumerate}
\end{prop}

\begin{proof}\sloppy (1) By Proposition \ref{prop:tpalltails}
there is some $u\geq 1$ such that $t((R^{2n})^u)=t((R^{2n})^{(u+1)})$,
and so $r_{2nu-1}=r_{2n(u+1)-1}$ where 
$(r_2,r_3,\dots,r_{2n(u+1)})$ is the maximal left overlap sequence for 
$(R^{2n})^{(u+1)}$. Since $R^{2n}$ is tightly packed, 
$r_{2n} = r_{2nu} = r_{2n(u+1)}$.  Thus
$R^{2n}$ is left $(2nu-1)$-stable.

Now suppose that $R^{2n}$ is left $j$-stable, and let $v$ be such that 
$j< 2nv$. Let $(r_2,r_3,\dots,r_{2n(v+1)})$ be the maximal left overlap 
sequence for $(R^{2n})^{(v+1)}$.  Then by Proposition \ref{conv:p1}, 
we have $r_k = r_{2n+k}$ for all $k$, with $j \leq k \leq 2nv$.
In particular, $r_{2nv-1}=r_{2n(v+1)-1}$.  Since $R^{2n}$ is tightly packed, 
$r_{2nv}=r_{2n(v+1)}$, and so $t((R^{2n})^v) = t((R^{2n})^{(v+1)})$.

The proof of (2) is similar and is left to the reader.
\end{proof}

\begin{defin} Let $R^{2n}$ be an extending closed walk.  We say 
$R^{2n}$ {\it stabilizes} at $u \geq 1$ if $R^{2n}$ is left and right
$j$-stable for some $j\leq 2nu-1$.
\end{defin}

\begin{prop}\label{conv:propremark} 
Let $R^{2n}$ be a tightly packed extending closed walk which stabilizes at 
$u \geq 1$. Then, for all $k\ge 1$ and $0 \leq i < 2n$,
\begin{enumerate}
\item if $(R^{2n})^{(u+k)} = R^{2n(u+k-1)+i}q$ and 
$(R^{2n})^{2u} = R^{2nu+i}q^{\prime}$ where
$R^{2n(u+k-1)+i}\in \R^{2n(u+k-1)+i}$ and $R^{2n+i}\in \R^{2n+i}$ then
$$t(R^{2n(u+k-1)+i})= t(R^{2nu+i}),$$ and 
\item if $(R^{2n})^{(u+k)} = qR^{2n(u+k-1)+i}$ and 
$(R^{2n})^{2u} = q^{\prime}R^{2nu+i}$ where
$R^{2n(u+k-1)+i}\in \R^{2n(u+k-1)+i}$ and $R^{2n+i}\in \R^{2n+i}$ then
$$b(R^{2n(u+k-1)+i})= b(R^{2nu+i}).$$  
\end{enumerate}
\end{prop}

\begin{proof}
By Proposition \ref{conv:p2}, $r_{2nu-1} =r_{2n(u+k)-1}$ 
and $r_{2nu}=r_{2n(u+k)}$ for all $k\ge 1$, where
$(r_2,r_3,\dots,r_{2n(u+k)})$ is the maximal left
overlap sequence for $(R^{2n})^{(u+k)}$.  From Proposition \ref{conv:p1},
we have $r_{2nu+i}=r_{2nu+2nk+i}$ for all $k\ge 1$ and 
$0 \leq i < 2n$. 
But then, with the above notation, and since 
$(r_2,r_3,\dots,r_{2n(u+k-1)+i})$ is a maximal left overlap sequence for 
$R^{2n(u+k-1)+i}$ and 
$(r_2,r_3,\dots,r_{2nu+i})$ is a maximal left overlap sequence for 
$R^{2nu+i}$, it follows that
$t(R^{2n(u+k-1)+i})= t(R^{2nu+i})$ for all $k\ge 1$ and 
$0 \leq i < 2n$.  

The proof of (2) for beginning paths is similar.
\end{proof}

\section{The equivalence relation}

In this section we begin our study of the graded centre of the Ext algebra 
modulo the ideal generated by the homogeneous nilpotent elements. 
The graded centre of the Ext algebra is denoted $Z_\gr(E(\L))$ and 
is the subring of $E(\L)$ generated by all homogeneous elements $z$ 
such that $zg = (-1)^{|z||g|}gz$ for each homogeneous element $g$ in $E(\L)$.

Let $Z^m(E(\L))$ denote the vector space of homogeneous elements of degree 
$m$ in the graded centre $Z_\gr(E(\L))$. Thus
$Z_\gr(E(\L)) = \oplus_{m\geq 0}Z^m(E(\L))$. 
For ease of notation we do not repeat the subscript $_\gr$ on the
homogeneous spaces $Z^m(E(\L))$; note that we do not use the notation 
$Z$ to represent any ungraded centre of a ring.

We consider here elements in $Z_\gr(E(\L))$ of even degree, noting that if the
characteristic of the field is not 2, then every element of $Z_\gr(E(\L))$
is of even degree. We return to this point with a discussion of elements of 
odd degree in $Z_\gr(E(\L))$ at the end of section 4.

Let $z$ be a homogeneous element of degree $2n$ in
$Z_\gr(E(\L))$, that is, let $z$ be in $Z^{2n}(E(\L))$. 
Then, using the multiplicative basis for
$E(\L)$, we may write $z = \sum_{i=1}^l\alpha_ig_i^{2n}$ for some
$\alpha_i \in K\setminus\{0\}$. For each $i = 1, \ldots , l$, we
say that $g_i^{2n}$ {\it occurs} in $z$.

The element $z$ is in $Z_\gr(E(\L))$ and so for each basis element
$g^0_j$ of $\G^0$, we have $g^0_jz=zg^0_j$. Thus, for each $i=1,
\ldots , l$, there is some vertex $e_{v_i}$ with
$R_i^{2n}=e_{v_i}R_i^{2n}e_{v_i}$. Hence we may write
$z=\sum_{i=1}^l \alpha_ig_i^{2n}$, where the $R_i^{2n}$ are closed
walks in the quiver, and $\alpha_i \in K\setminus\{0\}$.

The following result and the ideas in its proof are used throughout this
paper.

\begin{lem}\label{lemma:zero}
If $i \neq j$ then $g_i^{2n}g_j^{2n} = 0$ whenever $g_i^{2n}$ and
$g_j^{2n}$ occur in $z$.
\end{lem}

\begin{proof}
Suppose that $g_i^{2n}g_j^{2n} \neq 0$. We show first that
$zg_j^{2n}\neq 0$. For suppose for contradiction that $zg_j^{2n} =
0$. Then, since $z=\sum_{k=1}^l \alpha_kg_k^{2n}$ with $\alpha_k
\in K\setminus\{0\}$, we have
$\sum_{k=1}^l\alpha_kg_k^{2n}g_j^{2n} = 0$. Thus
$\sum_{g_k^{2n}g_j^{2n} \neq 0} \alpha_kg_k^{2n}g_j^{2n} = 0$, and
this is not the empty sum since $g_i^{2n}g_j^{2n} \neq 0$. Now
each of the terms $g_k^{2n}g_j^{2n}$ where $g_k^{2n}g_j^{2n} \neq
0$, is a basis element in $E(\L)$ since we have a multiplicative
basis. Thus $\{g_k^{2n}g_j^{2n}: g_k^{2n}g_j^{2n} \neq 0\}$ is a
linearly independent set by Proposition \ref{prop:extcancel}. 
Hence each coefficient $\alpha_k$ in the
sum $\sum_{g_k^{2n}g_j^{2n} \neq 0} \alpha_kg_k^{2n}g_j^{2n}$ must
be zero, and in particular $\alpha_i = 0$. But this is a
contradiction. Thus $zg_j^{2n}\neq 0$.

Now $z \in Z_\gr(E(\L))$, so $zg_j^{2n} = g_j^{2n}z$. Thus
$$\sum_{k=1}^l \alpha_kg_k^{2n}g_j^{2n} =
\sum_{k=1}^l \alpha_kg_j^{2n}g_k^{2n}$$
and so
$$\sum_{g_k^{2n}g_j^{2n} \neq 0} \alpha_kg_k^{2n}g_j^{2n} =
\sum_{g_j^{2n}g_k^{2n} \neq 0} \alpha_kg_j^{2n}g_k^{2n}.$$
Now each $g_k^{2n}g_j^{2n}$ and $g_j^{2n}g_k^{2n}$ in this last equality is
a basis element in $E(\L)$ so since $g_i^{2n}g_j^{2n} \neq 0$,
there is some unique $i'$ with $g_i^{2n}g_j^{2n} = g_j^{2n}g_{i'}^{2n}$
and $\alpha_i = \alpha_{i'}$.
Hence from Proposition \ref{prop:equal} we have
$R_i^{2n} = R_j^{2n}$ and so $i = j$.

Thus if $i \neq j$ then $g_i^{2n}g_j^{2n} = 0$.
\end{proof}

Now we introduce an equivalence relation on certain closed walks in
${\mathcal Q}$.

\begin{defin}\label{defin:liesonequiv}
\begin{enumerate}
\item Let $p$ be any path and let $q$ be a closed walk in a quiver ${\mathcal
Q}$. Then $p$ {\it lies on} $q$ if $p$ is a subpath of $q^s$ for
some $s \geq 1$.
\item Let $p$ and $q$ be closed walks in a quiver ${\mathcal
Q}$. We define a relation by
$p \sim q$ if $p$ lies on $q$ and $q$ lies on $p$.
\end{enumerate}
\end{defin}

\begin{remark}
(1) It is easy to verify that this is an equivalence relation.

(2) If $q$ and $q'$ are both in the same equivalence class, then a path $p$
lies on $q$ if and only if $p$ lies on $q'$.

(3) Throughout this paper, any reference to an equivalence relation is
referring to that of Definition \ref{defin:liesonequiv}.
\end{remark}

The next result uses the above definition to show that the problem of
describing the elements in  $Z_\gr(E(\L))$ can be reduced to each
equivalence class of closed walks. 

\begin{prop}\label{prop:closedwalks}
Let $z \in Z^{2n}(E(\L))$. Then
$z=\sum_{i=1}^r z_i$ with $z_i=\sum_{g_j^{2n}\in
\F_i}\alpha_jg_j^{2n}$, where $\alpha_j \in K\setminus\{0\}$, each
$\F_i$ is a subset of $\E_i$, and the sets $\E_i$ are distinct
equivalence classes of closed walks in $K{\mathcal Q}$. Moreover
each $z_i$ is in $Z_{\gr}(E(\L))$.
\end{prop}

\begin{proof}
Let $z \in Z^{2n}(E(\L))$. We have seen that we may write $z=\sum_{i=1}^l
\alpha_ig_i^{2n}$, where the $R_i^{2n}$ are closed walks in the
quiver and $\alpha_i \in K\setminus\{0\}$. Using the equivalence
relation in Definition \ref{defin:liesonequiv}, each $R_i^{2n}$ lies in
some equivalence class. Thus we may write $z=\sum_{i=1}^r z_i$
with $z_i=\sum_{g_j^{2n}\in \F_i}\alpha_jg_j^{2n}$, where each
$\F_i$ is a subset of $\E_i$, and the sets $\E_i$ are distinct
equivalence classes of closed walks.

Now we show that each $z_i$ is itself in $Z_\gr(E(\L))$. Let $g^m$
be an arbitrary basis element of $E(\L)$ of degree $m$ such that
$g^mz_i\neq 0$. Thus there exists $j$ so that $g_j^{2n}$ occurs in $z_i$ 
and $g^mg_j^{2n}\neq 0$.
Using the arguments in the proof of Lemma \ref{lemma:zero}, we
have that $g^mz \neq 0$, and so there exists $j'$ such that
$g_{j'}^{2n}$ occurs in $z$, $g^mg_j^{2n} = g_{j'}^{2n}g^m$ and
$\alpha_j = \alpha_{j'}$. Then, by Lemma \ref{lem:walk},
$R_j^{2n}$ and $R_{j'}^{2n}$ are in the same equivalence class
$\E_i$, and thus $g_{j'}^{2n}$ occurs in $z_i$. Continuing in this
way and using the corresponding argument with $z_ig^m \neq 0$
shows that nonzero terms in $g^mz_i$ and $z_ig^m$ pair up with the
same coefficients. Thus $g^mz_i = z_ig^m$. Hence $z_i \in
Z_\gr(E(\L))$ for all $i$.
\end{proof}

\begin{remark}
(1) Note that we do not assume that all $g_j^{2n}$, for which $R_j^{2n}$
lies in the
equivalence class $\E_i$, occur in $z$ with nonzero coefficient.

(2) It follows from Lemma \ref{lemma:zero} that $z_iz_j=0$ for
$i\neq j$ and thus $z^v = \sum_{i=1}^r z_i^v$ for any $v \geq 1$.
\end{remark}

From Proposition \ref{prop:closedwalks}, we may suppose without loss of
generality that if $z \in Z^{2n}(E(\L))$ then 
$z=\sum_{g_j^{2n}\in \F}\alpha_jg_j^{2n}$ where 
$\F$ is a subset of the equivalence class $\E$.
We say that such an element $z$ is 
{\it associated to the equivalence class} $\E$.

The next stage is to determine more precisely the subset $\F$ of the
equivalence class $\E$.

\section{Non-nilpotent elements of the graded centre}

Let $\N_Z$ denote the ideal of $Z_\gr(E(\L))$ generated by the
homogeneous nilpotent elements. Then $\N_Z$ is a graded ideal and we
write $\N_Z = \oplus_{m\geq 0}\N^m_Z$, with $\N^m_Z$ being the vector
space of nilpotent homogeneous elements of degree $m$ in
$Z_\gr(E(\L))$. Now the quotient $Z_\gr(E(\L))/\N_Z$ has the
decomposition
$$Z_\gr(E(\L))/\N_Z = \oplus_{m\geq 0}\left ( Z^m(E(\L))/\N^m_Z\right ).$$

If char$K \neq 2$, then all elements of odd degree in $Z_\gr(E(\L))$ are 
nilpotent. Thus, in every characteristic, $Z_\gr(E(\L))/\N_Z$ is a commutative
$K$-algebra.

Fix some equivalence class $\E$. Again we consider elements in
$Z_\gr(E(\L))$ of even degree. Using Proposition
\ref{prop:closedwalks}, it is enough to consider elements $z\in
Z^{2n}(E(\L))$ which are non-nilpotent and are associated to the
equivalence class $\E$. Write $z = \sum_{i=1}^l\alpha_ig_i^{2n}$ with
$\alpha_i \in K\setminus\{0\}$ and $R_i^{2n} \in \E$.

Now partition the set $\{1, \ldots , l\}$ into two subsets $\I$
and $\I '$ where $\I = \{i \colon (g_i^{2n})^k \neq 0 \mbox{ for
all } k \geq 1\}$ and $\I ' = \{i \colon (g_i^{2n})^k = 0 \mbox{
for some } k \geq 1\}$. Then we may write $z = \sum_{i \in \I}
\alpha_ig_i^{2n} + y$ where $y = \sum_{i \in \I '}
\alpha_ig_i^{2n}$. Using Lemma \ref{lemma:zero}, $y$ is nilpotent.

Now, from Proposition \ref{prop:tpalltails}, we may choose $u$ large enough so
that $y^u = 0$ and, for each $i \in \I$, we have
$t((R_i^{2n})^u) = t((R_i^{2n})^{(u+1)})$ and 
$b((R_i^{2n})^u) = b((R_i^{2n})^{(u+1)})$.
Then $z^u = \sum_{i \in \I} \alpha_i^ug_i^{2nu}$ (from Lemma
\ref{lemma:zero}). Note that $z^u$ is also in $Z_\gr(E(\L))$. We keep this
notation for the rest of this section.

The following example shows an element of the graded centre for which we
require $u=2$.

\begin{example}\label{ex:gradedcentre}
Let $\L = K{\mathcal Q}/I$ where ${\mathcal Q}$ is the quiver
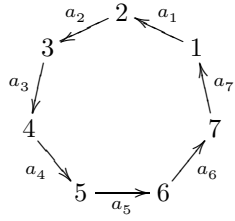
\[
\def\alphanum{\ifcase\xypolynode\or 1\or 2
\or 3\or 4\or 5\or 6\or 7\fi}
\xy/r3pc/: {\xypolygon7{~<<{@{}}~><{}
~*{\alphanum}
~>>{_{a_\xypolynode^{}}}}}
\endxy\]
and 
$$I = \langle a_1a_2a_3a_4, a_2a_3a_4a_5, a_4a_5a_6, a_5a_6a_7,
a_6a_7a_1a_2, a_7a_1a_2a_3\rangle.$$

\bigskip

The elements of the set $\R^4$ are
$$\begin{array}{lcl}
R_1^4 = a_1a_2a_3a_4a_5a_6a_7 & \mbox{ with } & t(R_1^4) = a_6a_7\\
R_2^4 = a_4a_5a_6a_7a_1a_2a_3 & \mbox{ with } & t(R_2^4) = a_1a_2a_3\\
R_3^4 = a_5a_6a_7a_1a_2a_3a_4 & \mbox{ with } & t(R_3^4) = a_3a_4\\
R_4^4 = a_7a_1a_2a_3a_4a_5a_6 & \mbox{ with } & t(R_4^4) = a_5a_6\\
R_5^4 = a_2a_3a_4a_5a_6a_7a_1a_2 & \mbox{ with } & t(R_5^4) = a_7a_1a_2.
\end{array}$$
It may be verified that $R_1^4, R_2^4, R_3^4$ and $R_4^4$ are all tightly
packed extending closed walks. Note that $R_5^4$ is not a closed path so 
the element $g_5^4$ is nilpotent. Let $z = g_1^4+g_2^4+g_3^4+g_4^4$. Then $z
\in Z_\gr(E(\L))$ and is non-nilpotent.

Now
$$\R^8 = \{ (R_1^4)^2, (R_2^4)^2, (R_3^4)^2, (R_4^4)^2\}$$
with
$$t((R_1^4)^2) = a_7, t((R_2^4)^2) = a_1a_2a_3, t((R_3^4)^2) = a_4, 
t((R_1^4)^2) = a_5a_6.$$
Thus if $u = 1$, then $t((R_1^4)^u) \neq t((R_1^4)^{(u+1)})$. However, for 
$s \geq 3$, if $(r_2, r_3, \ldots , r_{4s})$ is a maximal left overlap
sequence for $(R_1^4)^s$ then $r_7 = r_{11} = a_4a_5a_6$ and $r_8 = r_{12} =
a_5a_6a_7$ and so $R_1^4$ is left 7-stable. Moreover $R_1^4$ is also right 
7-stable. Indeed, each of 
$R_1^4, R_2^4, R_3^4$ and $R_4^4$ is both left 7-stable and right 7-stable.
So, for $u = 2$, we have $t((R_i^4)^u) = t((R_i^4)^{(u+1)})$ for $i = 1,
\ldots , 4$.

Note that in this example $Z_\gr(E(\L))/\N_Z \cong K[\overline{z}]$.
\end{example}

\begin{remark}
Suppose that $g_i^{2n}$ is in $\G^{2n}$ and that 
$(g_i^{2n})^u \neq 0$ for some $u \geq 1$. Then $(g_i^{2n})^u$ is in $\G^{2nu}$
and so $(R_i^{2n})^u$ is the path 
corresponding to some element of $\R^{2nu}$. We write $R_i^{2nu}$ for the
element of $\R^{2nu}$ with path $(R_i^{2n})^u$, and will use this notation 
throughout the rest of the paper.
\end{remark}

The next results describe the paths $R_i^{2nu}$ for $i \in \I$.

\begin{prop}\label{prop:reltightlypacked}
Suppose $g_i^{2nu}$ occurs in $z^u$. Then the following statements
hold.
\begin{enumerate}
\item The path $R_i^{2nu}$ is tightly packed and may be expressed in the 
form
$$R_i^{2nu} = r_{i,2}r_{i,4}\cdots r_{i,2nu}$$
where $r_{i,2l}$ is the path corresponding to some basis element
$g_{i,2l}^2 \in \G^2$. Thus the $r_{i,2l}$ are relations lying on
every representative of $\E$.
\item $R_i^{2n}$ is tightly packed with
$$R_i^{2n} = r_{i,2}r_{i,4}\cdots r_{i,2n}.$$
\item The relations in \emph{(1)} are such that $r_{i,2cn+2d} = r_{i,2d}$ for
$c = 1, \ldots, u-1, d = 1, \ldots, n-1$.
\end{enumerate}
\end{prop}

\begin{proof}
(1) Without loss of generality suppose $i = 1$. Let $r_{1,2}$ be
the first relation in $R_1^{2nu}$ and let $g_{1,2}^2$ be the basis
element in $E(\L)$ with path $r_{1,2}$. Then $(g_1^{2nu})^2 \neq
0$ and so Proposition \ref{prop:multrelation} gives $g_1^{2nu}g_{1,2}^2 \neq
0$. Thus, using the ideas in the proof of Lemma \ref{lemma:zero},
$z^ug_{1,2}^2 \neq 0$ and there is some $i\in\I$ with $g_1^{2nu}g_{1,2}^2 =
g_{1,2}^2g_i^{2nu}$. Hence $R_1^{2nu}r_{1,2} = r_{1,2}R_i^{2nu}$.
Now the first relation of $R_i^{2nu}$ is $r_{i,2}$, and thus
$r_{1,2}R_i^{2nu}$ is a path with prefix $r_{1,2}r_{i,2}$. Hence
$R_1^{2nu}$ has prefix $r_{1,2}r_{i,2}$. Let $r_{1,4} = r_{i,2}$ so that
$R_1^{2nu}$ has prefix $r_{1,2}r_{1,4}$.

Continuing inductively
shows that $R_1^{2nu} = r_{1,2} \cdots r_{1,2nu}p$ for some path $p$, and
that there is a maximal left overlap sequence with path 
$r_{1,2}r_{1,4} \cdots r_{1,2nu}$. By Proposition \ref{prop:tpleftoverlap}, 
$r_{1,2}r_{1,4} \cdots r_{1,2nu} \in \R^{2nu}$ and so, by
Proposition \ref{prop:uniquepaths}, the path $p$ is a vertex. Thus 
$R_1^{2nu} = r_{1,2} \cdots r_{1,2nu}$ and so
$R_1^{2nu}$ is a product of relations. Hence, for each
$g_i^{2nu}$ occurring in $z^u$, $R_i^{2nu}$ is a product of
relations and we write $R_i^{2nu} = r_{i,2}r_{i,4}\cdots r_{i,
2nu}$. Thus $R_i^{2nu}$ is tightly packed.

(2) and (3) now follow from Proposition \ref{prop:tppowers}.
\end{proof}

\begin{lem}\label{lem:relationcoeff}
With the notation of \emph{Proposition \ref{prop:reltightlypacked}},
let $r_{i,2h}$ be a relation in $R_i^{2n}$, for some
$i\in\I$ and $h=1,\ldots , n$.
Then $r_{i,2h}$ is the first relation of some $R_j^{2n}$ with $j\in\I$.
Moreover $\alpha_i = \alpha_j$.
\end{lem}

\begin{proof}
Since $i\in\I$, the element $g_i^{2nu}$ occurs in $z^u$ and
$g_i^{2nu}g_i^{2n} \neq 0$. Thus from Proposition \ref{prop:multrelation},
$g_i^{2nu}g^2_{i,2} \neq 0$ where $g^2_{i,2}$ is the basis element
corresponding to the relation $r_{i,2}$. Using the arguments of Lemma
\ref{lemma:zero}, we have that $z^ug^2_{i,2} \neq 0$ and that there is some
$j \in \I$ with $g_i^{2nu}g^2_{i,2} = g^2_{i,2}g_j^{2nu}$ and $\alpha_i =
\alpha_j$. From Proposition \ref{prop:reltightlypacked}, $R_i^{2nu} =
(r_{i,2}r_{i,4}\cdots r_{i,2n})^u$ and so
$R_j^{2n} = r_{i,4}r_{i,6}\cdots r_{i, 2n}r_{i,2}$. Thus $r_{i,4}$ is the
first relation of some $R_j^{2n}$ with $j\in\I$ and $\alpha_i = \alpha_j$.
The result now follows by induction.
\end{proof}

\begin{remark}
In the above proof, we have also shown that $r_{i,2}$ is the last relation 
of $R_j^{2n}$ with $j\in\I$ and $\alpha_i = \alpha_j$. Thus by induction, 
each $r_{i,2h}$ is also the last relation of some $R_j^{2n}$ with $j\in\I$.
\end{remark}

\begin{lem}\label{lem:firstlastrelation}
With the notation of \emph{Proposition \ref{prop:reltightlypacked}},
let $R_i^{2n} = r_{i,2}r_{i,4}\cdots r_{i, 2n}$ with $i\in\I$. If $r \in \R^2$
with $R_i^{2n}r$ in $\R^{2n+2}$, then $r = r_{i,2}$.
\end{lem}

\begin{proof}
Let $g^2$ be the basis element in $E(\L)$ with path $r$. Then $g_i^{2n}g^2
\neq 0$.
With the arguments of Lemma \ref{lemma:zero}, we have that $zg^2 \neq 0$ 
and that there is some $j \in \{ 1, \ldots , l\}$ with $g_i^{2n}g^2 = 
g^2g_j^{2n}$. Hence $r$ is the first relation in $R_i^{2n}$, that is, 
$r = r_{i,2}$.
\end{proof}

\begin{thm}\label{thm:commonrelation}\sloppy 
With the notation of \emph{Proposition \ref{prop:reltightlypacked}}, 
suppose that $R_1^{2n} = r_2r_4\cdots r_{2n}$, 
$R_2^{2n} = \tilde{r}_2\tilde{r}_4\cdots \tilde{r}_{2n}$ and that 
$r_{2n} = \tilde{r}_{2n}$. Then $R_1^{2n} = R_2^{2n}$.
\end{thm}

\begin{proof}
We begin by showing that $r_2 =\tilde{r}_2$.
The tails $t(R_1^{2n})$ and $t(R_2^{2n})$ are both suffixes of 
$r_{2n}=\tilde r_{2n}$, so without loss of generality, we may assume 
that the length of $t(R_1^{2n})$ is
greater than or equal to the length of $t(R_2^{2n})$.
Let $(r_2,r_3,\ldots, r_{2n})$ be the maximal left 
overlap sequence for $R_1^{2n}$ and let 
$(\tilde{r}_2,\tilde{r}_3,\ldots, \tilde{r}_{4n})$ be the maximal left 
overlap sequence for $(R_2^{2n})^2$.  
Note that $\tilde{r}_{2n+2i} = \tilde{r}_{2i}$ for $i = 1,\ldots , n$.
Now $\tilde{r}_{2n+1}$ overlaps $t(R^{2n}_2)$ and so
$\tilde{r}_{2n+1}$ overlaps $t(R^{2n}_1)$.  
Since $r_{2n} = \tilde{r}_{2n}$, we have that 
$(r_2,r_3,\ldots,r_{2n},\tilde{r}_{2n+1},\tilde{r}_{2})$ is a left overlap
sequence with path $r_2\cdots r_{2n}\tilde{r}_{2}$. Since this path is 
tightly packed, it follows from Proposition \ref{prop:tpleftoverlap} 
that $R^{2n}_1\tilde{r}_{2} \in \R^{2n+2}$.  
Hence from Lemma \ref{lem:firstlastrelation} we have $\tilde{r}_2 = r_2$.  

Now from the proof of Lemma \ref{lem:relationcoeff} and the remark
following, we may repeat this process with the terms $r_4r_6\cdots r_{2n}r_2$ 
and $\tilde{r}_4\tilde{r}_6\cdots \tilde{r}_{2n}\tilde{r}_2$ where 
$\tilde{r}_{2n} = r_{2n}$ and $\tilde{r}_2 = r_2$ to show that 
$\tilde{r}_4 = r_4$. Hence by induction we have $R^{2n}_1 = R^{2n}_2$.  
\end{proof}

We now show that the coefficients $\alpha_i$ for $i \in \I$ are
all equal.

\begin{prop}\label{prop:coeff}
Let $z \in Z^{2n}(E(\L))$ be a non-nilpotent element associated to the
equivalence class $\E$.
Write $z = \sum_{i\in\I} \alpha_ig_i^{2n} + y$ where y is nilpotent,
$\alpha_i \in K\setminus\{0\}$ and all $R_i^{2n}$ lie in $\E$. 
Then all the coefficients $\alpha_i$ for $i \in \I$ are
equal, so $z = \alpha\sum_{i\in\I} g_i^{2n} + y$ for some 
$\alpha \in K\setminus\{0\}$ where y is nilpotent. Moreover each path
$R_i^{2n}$ with $i \in \I$ has the same length.
\end{prop}

\begin{proof}
With appropriate relabeling of the set $\G^{2n}$, let $1 \in \I$,
and write $R_1^{2nu} = R_1^{2nu-1}t(R_1^{2nu})$.
Let $t_1 = t(R_1^{2nu})$. By definition
of the set $\I$, we have that $(g_1^{2n})^kg_1^{2nu}
\neq 0$ so that $R_1^{2nk}R_1^{2nu} = R_1^{2nk}R_1^{2nu-1}t_1$
for all $k\geq 1$.

We show first that $(g_1^{2n})^kg_1^{2nu-1}\neq 0$ for each $k
\geq 1$. Now $(g_1^{2n})^kg_1^{2nu} = g_1^{2n(k+u)}$ so that
$$\begin{array}{rcl}
R_1^{2n(k+u)} & = & R_1^{2n(k+u)-1}t(R_1^{2n(k+u)})\\
& = & R_1^{2n(k+u)-1}t(R_1^{2nu})\\
& = & R_1^{2n(k+u)-1}t_1
\end{array}$$
with the penultimate equality coming from the definition of $u$.
Thus $R_1^{2n(k+u)-1}t_1 = R_1^{2n(k+u)} =
R_1^{2nk}R_1^{2nu-1}t_1$. Hence $R_1^{2n(k+u)-1} =
R_1^{2nk}R_1^{2nu-1}$ so that
$$R_1^{2nk}R_1^{2nu-1} \in \R^{2n(u+k)-1}.$$
Thus $(g_1^{2n})^kg_1^{2nu-1}\neq 0$ for each $k\geq 1$.

Now with $k = 1$, $z \in Z_\gr(E(\L))$ so, using the arguments of
Lemma \ref{lemma:zero}, we have $zg_1^{2nu-1} \neq 0$ and, again relabeling
$\G^{2n}$ if necessary, $g_1^{2n}g_1^{2nu-1} =
g_1^{2nu-1}g_2^{2n}$ and $\alpha_1 = \alpha_2$. So $g_1^{2n}$
and $g_2^{2n}$ occur in $z$ with the same coefficients and
$\ell(R_1^{2n}) = \ell(R_2^{2n})$.

We show next that $2 \in \I$, that is, $g_2^{2n}$ is not nilpotent.
For each $k \geq 1$, $z^k \in
Z_\gr(E(\L))$ so, using the arguments of Lemma \ref{lemma:zero},
we have $z^kg_1^{2nu-1} \neq 0$ and there is some $j(k)$ with
$(g_1^{2n})^kg_1^{2nu-1} = g_1^{2nu-1}(g_{j(k)}^{2n})^k$.
But $g_1^{2n}g_1^{2nu-1} = g_1^{2nu-1}g_2^{2n}$, so by
uniqueness of paths it follows that $j(k) = 2$ for all $k$. Hence
$$(g_1^{2n})^kg_1^{2nu-1} = g_1^{2nu-1}(g_2^{2n})^k
\mbox{ for all $k \geq 1$.}$$
Thus
$g_1^{2nu-1}(g_2^{2n})^k\neq 0$ and hence $(g_2^{2n})^k \neq
0$ for all $k \geq 1$. Thus $2 \in \I$.

Now we show that in fact $R_2^{2nu} = t_1R_1^{2nu-1}$. With $k
= u$, $(g_1^{2n})^ug_1^{2nu-1} = g_1^{2nu-1}(g_2^{2n})^u$ so that
$$R_1^{2nu}R_1^{2nu-1} = R_1^{2nu-1}R_2^{2nu}.$$
Write $R_2^{2nu} = qR_2^{2nu-1}$ for some path $q$. Then we have
$$R_1^{2nu}R_1^{2nu-1} = R_1^{2nu-1}qR_2^{2nu-1}.$$
By uniqueness of paths, $R_1^{2nu-1} = R_2^{2nu-1}$ and
$R_1^{2nu} = R_1^{2nu-1}q$.
Since $R_1^{2nu} = R_1^{2nu-1}t_1$ it follows that $t_1 = q$.
Thus $R_2^{2nu} = t_1R_1^{2nu-1}$.

So, proceeding inductively, we have
$$\begin{array}{rclcl}
R_1^{2nu} & = & R_1^{2nu-1}t_1 & & \\
R_2^{2nu} & = & t_1R_1^{2nu-1} & = & R_2^{2nu-1}t_2\\
R_3^{2nu} & = & t_2R_2^{2nu-1} & = & R_3^{2nu-1}t_3\\
& \vdots & & \\
R_f^{2nu} & = & t_{f-1}R_{f-1}^{2nu-1} & = & R_f^{2nu-1}t_f\\
& \vdots & &
\end{array}$$
with $\{1, \ldots , f\} \subseteq \I$. Moreover $\alpha_1 =
\alpha_2 = \cdots = \alpha_f$ and $\ell(R_1^{2n}) = \ell(R_2^{2n})
= \cdots = \ell(R_f^{2n})$.

Since the set $I$ is finite, there exists a minimal $f$ where this
list of $R_i^{2nu}$ repeats, that is, there exists a minimal $f$
such that $R_{f+1}^{2nu} = R_s^{2nu}$ for some $1 \leq s \leq f$.
We show that $s = 1$. Suppose for contradiction that $s > 1$. Then
$t_fR_f^{2nu-1} = R_{f+1}^{2nu} = R_s^{2nu} =
t_{s-1}R_{s-1}^{2nu-1}$. So by Proposition \ref{prop:uniquepaths}, it
follows that $R_f^{2nu-1} = R_{s-1}^{2nu-1}$ and thus $t_f =
t_{s-1}$. Hence $R_f^{2nu} = R_f^{2nu-1}t_f =
R_{s-1}^{2nu-1}t_{s-1} = R_{s-1}^{2nu}$, which contradicts the
minimality of $s$. Thus $s = 1$, and so $R_{f+1}^{2nu} =
R_1^{2nu}$.

We have that $\alpha_1 = \cdots = \alpha_f$ so write $\alpha = \alpha_1$.
Then $g_1^{2n}, \ldots , g_f^{2n}$ all occur in $z$ with
coefficient $\alpha$. Moreover from Lemma \ref{lem:relationcoeff},
for each relation ${r_{i,2h}}$ in some $R_1^{2n}, \ldots , R_f^{2n}$
there is some $R_j^{2n}$ with $j\in\I$ where $g_j^{2n}$ occurs in $z$
with the same coefficient $\alpha$. We now let $\J_1$ denote the subset
of $\I$ consisting of all $i$ such that
$g_i^{2n}$ occurs in $z$ with coefficient $\alpha$. Continuing in this way,
we partition the set $\I$ into disjoint
sets $\J_1, \ldots , \J_d$ with respect to this construction.

Then $z = \sum_{j=1}^d\alpha_j(\sum_{i \in \J_j}g_i^{2n}) + y$, where $y$ is
nilpotent. Let $w_j = \sum_{i \in \J_j}g_i^{2n}$ so that $z =
\sum_{j=1}^d\alpha_jw_j + y$ with $y$ nilpotent and $\alpha_j \in
K\setminus\{0\}$.

\sloppy We complete the proof once we have shown that $d = 1$, since
then $z = \alpha\sum_{i\in\I} g_i^{2n} + y$, where $y$ is nilpotent
and $\alpha \in K\setminus\{0\}$. To do this we require the following
two propositions; the proof that $d = 1$ is the final part of the
second proposition.
\end{proof}

The next result describes how the different $R^{2n}$ are related for
$g^{2n}$ occurring in $z$.

\begin{prop}\label{relationsinmiddle}
Let $i\in\I$. With the notation of \emph{Proposition
\ref{prop:reltightlypacked}}
$$R_i^{2nu} = a_ir_{i+1,4}r_{i+1,6}\cdots r_{i+1, 2nu}t_i$$ and
$$R_i^{2nu} = t_{i-1}r_{i-1,2}r_{i-1,4}\cdots
r_{i-1, 2nu-2}\tilde{a}_{i-1}$$ for paths $a_i$ and $\tilde{a}_i$ with
$t_ia_i = r_{i+1, 2}$ and $\tilde{a}_{i-1}t_{i-1} = r_{i-1, 2nu}$.
\end{prop}

\begin{proof}
Without loss of generality, consider the case $i=1$.

Recall that $R_2^{2nu} = t_1R_1^{2nu-1}$ and $R_2^{2nu} =
r_{2,2}r_{2,4}\cdots r_{2, 2nu}$ so that $R_1^{2n-1}=
a_1r_{2,4}r_{2,6}\cdots r_{2, 2nu}$ for some path $a_1$ with
$t_1a_1 = r_{2, 2}$. Now $R_1^{2nu} = R_1^{2nu-1}t_1$. Thus
$R_1^{2nu} = a_1r_{2,4}r_{2,6}\cdots r_{2, 2nu}t_1$.

Similarly using $R_f^{2nu} = R_f^{2nu-1}t_f$ and $R_f^{2nu} =
r_{f,2}r_{f,4}\cdots r_{f, 2nu}$ we have $R_f^{2nu-1}=
r_{f,2}r_{f,4}\cdots r_{f, 2nu-2}\tilde{a}_f$ for some path
$\tilde{a}_f$ with $\tilde{a}_ft_f = r_{f,2nu}$. Now $R_1^{2nu} =
R_{f+1}^{2nu} = t_fR_f^{2nu-1}$. Thus $R_1^{2nu} =
t_fr_{f,2}r_{f,4}\cdots r_{f, 2nu-2}\tilde{a}_f$.
\end{proof}

The next result is the last part of the proof of Proposition
\ref{prop:coeff}.

\begin{prop}\label{prop:tpsequences}
With the notation already introduced, let $p=s_2s_4\cdots s_{2n}$ be a 
tightly packed closed walk in the equivalence class $\E$, with $s_2,
s_4, \ldots , s_{2n} \in\R^2$.
Then there is some $j\in\J_1$ such that $p = R_j^{2n}$ and, in particular,  
$p$ is in $\R^{2n}$. Hence $d = 1$.
\end{prop}

\begin{proof}
With appropriate relabeling of the set $\G^{2n}$, suppose that $1\in\J_1$ 
and write $R^{2n}_1 = r_{1,2}\dots r_{1,2n}$.  Consider first all tightly 
packed closed walks $p=s_2s_4\cdots s_{2n}$ which are in the 
equivalence class $\E$ and whose first relation $s_2$ 
either is $r_{1,2}$ or overlaps $r_{1,2}$.  Label these
paths $R^{2n}[1], R^{2n}[2],\ldots,R^{2n}[k]$ such that
\begin{enumerate}
\item $R^{2n}[1]= R^{2n}_1$, and
\item \sloppy $r[i+1]_2$ overlaps  $r[i]_2$, for $i=1,\dots,k-1$, where
we write $R^{2n}[i] = r[i]_2r[i]_4 \cdots r[i]_{2n}$.
\end{enumerate}
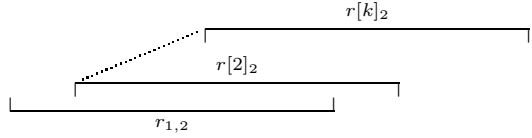
\[\xymatrix@W=0pt@M=0.3pt{%
\ar@{^{|}-^{|}}@<-1.25ex>[rrrrr]_{r_{1,2}} &
\ar@{_{|}-_{|}}@<1.25ex>[rrrrr]^{r[2]_2}\ar@{}@<1.25ex>[rrrrr]^(0){}="a"
& &
\ar@{_{|}-_{|}}@<6ex>[rrrrr]^{r[k]_2}\ar@{}@<6ex>[rrrrr]^(0){}="b"
&&&&&
\ar@{{}..{}} "a";"b" }\] %

Let $(r_{1,2},r_{1,3},\ldots,r_{1,2nu})$ be the maximal left overlap sequence
for $(R^{2n}_1)^u$. Then there is a left overlap sequence 
$(r[i]_2, r_{1,4}, r[i]_4, r_{1,6}, \ldots , r[i]_{2n})$ 
with path $R^{2n}[i]$ since $R^{2n}_1$ is tightly packed. So by 
Proposition \ref{prop:tpleftoverlap}, the path $R^{2n}[i]$ is in $\R^{2n}$
for $i=1,\ldots , k$.  

We now show that $R^{2n}[2] = R^{2n}_j$ for some $j\in \J_1$. 
To ease notation, let $r[2]_{2i} = \tilde{r}_{2i}$ for $i = 1, \ldots , n$.
If $\tilde{r}_2 = r_{1,2}$ then $R^{2n}[2] = R^{2n}_1$ and we are done.
Thus suppose $\tilde{r}_2 \neq r_{1,2}$. There are two cases to consider.  

First we suppose that $\tilde{r}_{2i} =  r_{1,2nl+2i+1}$ for some $i, l$ with
$1\le i\le n-1$ and $1\le l\le u-1$. Since $R_1^{2n}$ is tightly
packed, the relation $\tilde{r}_{2i}$ is situated as follows.
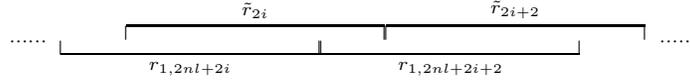
\[\xymatrix@W=0pt@M=0.3pt{%
\ar@{}@<-1.25ex>[r]^{\cdots\cdots} &
\ar@{^{|}-^{|}}@<-1.25ex>[rrrr]_{r_{1,2nl+2i}} &
\ar@{_{|}-_{|}}@<1.25ex>[rrrr]^{\tilde{r}_{2i}} &&&
\ar@{^{|}-^{|}}@<-1.25ex>[rrrr]_{r_{1,2nl+2i+2}} &
\ar@{_{|}-_{|}}@<1.25ex>[rrrr]^{\tilde{r}_{2i+2}} & & & &
\ar@{}@<-1.25ex>[r]^{\cdots\cdots} & }\] %

Now $R^{2n}[2]$ is also tightly
packed, so that $\tilde{r}_{2i+2} =  r_{1,2nl+2i+3}$.
Continuing in this way we have that the maximal left overlap sequence for
$(R_1^{2n})^u = R_1^{2nu}$ has the form
$$(r_{1,2},r_{1,3},\ldots,r_{1,2nl+2i}, \tilde{r}_{2i}, r_{1,2nl+2i+2},
\tilde{r}_{2i+2},\dots, \tilde{r}_{2n(u-l)-2}, r_{1,2nu}).$$
Now, 
$$(r_{1,2},r_{1,3},\ldots,r_{1,2nl+2i}, \tilde{r}_{2i}, r_{1,2nl+2i+2},
\tilde{r}_{2i+2},\dots, \tilde{r}_{2n(u-l)-2})$$
is a maximal left overlap sequence for $R_1^{2nu-1}$.  
By choice of $u$, $R_1^{2nu} = R_1^{2nu-1}t_1$, and so $t_1$ is the path
indicated below.
\[\xymatrix@W=0pt@M=0.3pt{%
\ar@{}@<-1.25ex>[r]^{\cdots\cdots} &
\ar@{^{|}-^{|}}@<-1.25ex>[rrrr]_{r_{1,2nu-2}} &
\ar@{_{|}-_{|}}@<1.25ex>[rrrr]^{\tilde{r}_{2n(u-l)-2}} &&&
\ar@{^{|}-^{|}}@<-1.25ex>[rrrr]_{r_{1,2nu}} &
\ar@{{<}-{>}}[rrr]^{t_1} & & & \ar@{}@<-1.25ex>[r]^{\cdots\cdots} & }\] %

{F}rom Proposition \ref{relationsinmiddle}, we have 
$$R_1^{2nu} = a_1r_{2,4}r_{2,6}\cdots r_{2, 2nu}t_1$$
and so $\tilde{r}_{2n(u-l)-2} = r_{2,2nu}$. Moreover, by definition of the
set $\J_1$ we have that $2 \in \J_1$.
Thus $r_{2,2} = r_{2,2nu+2} = \tilde{r}_{2n(u-l)}$.
Now $R^{2n}[2]$ is tightly packed and so $\tilde{r}_{2n(u-l)} =  \tilde{r}_{2n}$.
Hence $$R_2^{2nu} = \tilde{r}_{2n}\tilde{r}_{2}\cdots \tilde{r}_{2n-2}.$$
Using Lemma \ref{lem:relationcoeff}, it follows that
$\tilde{r}_{2}$ is the first relation of some $R^{2n}_j$ with $j\in\I$, that is,
$R_j^{2n} = \tilde{r}_{2}\tilde{r}_{4}\cdots \tilde{r}_{2n}$.
Again by definition of $\J_1$, we have that $j \in \J_1$.
Thus we have shown that $R^{2n}[2] = R^{2n}_j$
for some $j \in \J_1$.

For the second case, we suppose that $\tilde{r}_{2i} \neq  
r_{1,2nl+2i+1}$ for all $1\le i\le n-1$ and $1\le l\le u-1$. 

From our initial hypothesis, $\tilde{r}_{2}$ overlaps $r_{1,2}$. Now $r_{1,3}$ is
the first relation that overlaps $r_{1,2}$ since $(r_{1,2}, r_{1,3})$ 
is a maximal left overlap sequence. Thus $\tilde{r}_{2}$, which is the first
relation which starts a tightly packed closed walk and which overlaps
$r_{1,2}$, must overlap $r_{1,3}$. This is illustrated in the following diagram.
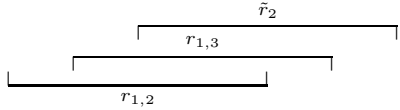
\[\xymatrix@W=0pt@M=0.3pt{%
\ar@{^{|}-^{|}}@<-1.25ex>[rrrr]_{r_{1,2}} &
\ar@{_{|}-_{|}}@<1.25ex>[rrrr]^{r_{1,3}} &
\ar@{_{|}-_{|}}@<4ex>[rrrr]^{\tilde{r}_2} & & & & }\] %

By maximality of the left overlap sequence $(r_2, r_3, \ldots , r_{2nu})$ for
$R_1^{2nu}$, the relation $r_{1,5}$ must begin at or before $\tilde{r}_4$. But by
the assumptions of this case, $r_{1,5} \neq \tilde{r}_4$ and so $\tilde{r}_4$ overlaps 
$r_{1,5}$ and starts strictly after $\mo(r_{1,5})$. 
Continuing in this way, we see that $\tilde{r}_{2k}$ overlaps
$r_{1,2k+1}$ for $k = 1 \ldots , n$.  Since $R^{2n}[2]$ is tightly packed,
$\tilde{r}_{2nl+2k} = \tilde{r}_{2k}$ for $k = 1 \ldots , n$, $l = 1 \ldots , u-1$.
Hence $\tilde{r}_{2k}$ overlaps $r_{1,2nl+2k+1}$ for $k = 1 \ldots , n$ and 
$l = 1 \ldots , u-1$.

Now consider the maximal left overlap sequence
$(r_{1,2},r_{1,3},\dots,r_{1,2n(u+1)})$ for $(R^{2n}_1)^{(u+1)}$.
As above, we have from Proposition \ref{relationsinmiddle}, that 
$$R_1^{2nu} = a_1r_{2,4}r_{2,6}\cdots r_{2, 2nu}t_1$$
and so $r_{1, 2nu-1} = r_{2,2nu}$. Thus
$$(r_{1,2nu-1},r_{1,2nu+1},\dots,r_{1,2nu+2n-3})$$ is a tightly packed path,
and hence
$$(r_{1,2nu+3},r_{1,2nu+3},\dots,r_{1,2nu+2n+1})$$ is also tightly packed.
Thus, since $R^{2n}[2]$ is tightly packed and $\tilde{r}_{2n-2}$ overlaps and 
is not equal to $r_{1,2nu-1}$, it follows that $\tilde{r}_{2n}$ overlaps and is 
not equal to $r_{1,2nu+1}$. So $\tilde{r}_{2k}$ overlaps and is not equal to 
$r_{1,2nu+2k+1}$ for $k = 1, \ldots , n$. Thus, recalling that
$r_{1,2nu+2k} = r_{1,2k}$ for $k = 1 \ldots , n$,
we have the following diagram.
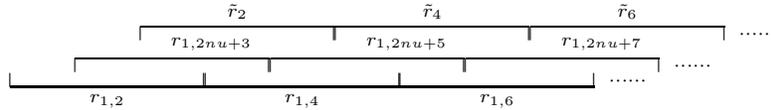
\[\xymatrix@W=0pt@M=0.3pt{%
\ar@{^{|}-^{|}}@<-1.25ex>[rrr]_{r_{1,2}} &
\ar@{_{|}-_{|}}@<1.25ex>[rrr]^(0.7){r_{1,2nu+3}} &
\ar@{_{|}-_{|}}@<4ex>[rrr]^{\tilde{r}_2} &
\ar@{^{|}-^{|}}@<-1.25ex>[rrr]_{r_{1,4}} &
\ar@{_{|}-_{|}}@<1.25ex>[rrr]^(0.7){r_{1,2nu+5}} &
\ar@{_{|}-_{|}}@<4ex>[rrr]^{\tilde{r}_4} &
\ar@{^{|}-^{|}}@<-1.25ex>[rrr]_{r_{1,6}} &
\ar@{_{|}-_{|}}@<1.25ex>[rrr]^(0.7){r_{1,2nu+7}} &
\ar@{_{|}-_{|}}@<4ex>[rrr]^{\tilde{r}_6} &
\ar@{}@<-1.8ex>[r]^{\cdots\cdots}
& \ar@{}@<-.5ex>[r]^{\cdots\cdots} & \ar@{}@<2.3ex>[r]^{\cdots\cdots} & }\] %
So we have distinct tightly packed sequences
$$(\tilde{r}_2,\tilde{r}_4,\dots,\tilde{r}_{2n}),$$ 
$$(r_{1,2nu+3},r_{1,2nu+5},\dots,r_{1,2nu+2n+1}),$$
and 
$$(r_{1,2},r_{1,4},\dots,r_{1,2n}),$$
for which $\tilde{r}_2$ overlaps $r_{1,2nu+3}$ and $r_{1,2nu+3}$ overlaps $r_{1,2}$.
But this contradicts (2) in the choice of labeling of $R^{2n}[2]$, and so 
this case cannot occur.

Thus we have shown that $R^{2n}[2] = R^{2n}_j$
for some $j\in \J_1$.   
It follows by induction, for each $i = 1, \ldots , k$, that
$R^{2n}[i] = R^{2n}_j$ for some $j\in \J_1$. 

Now let $p=s_2s_4\cdots s_{2n}$ be a tightly packed
path in the equivalence class $\E$.  Then there is some $1\le i\le n$, 
with either $s_2 = r_{1,2i}$ or $s_2$ overlaps $r_{1,2i}$.
By Lemma \ref{lem:relationcoeff}, $r_{1,2i}$ is the first relation of some
$R_c^{2n}$ with $c\in\J_1$. We now repeat the above arguments, replacing
$R_1^{2n}$ by $R_c^{2n}$, to show that there is some $j\in\J_1$ with
$p = R^{2n}_j$. Thus, in particular, $p \in \R^{2n}$.

Finally, suppose for contradiction that $d\geq 2$. Recalling the notation of
Proposition \ref{prop:coeff}, consider $w_1 = \sum_{i
\in \J_1}g_i^{2n}$ and $w_2 = \sum_{i \in \J_2}g_i^{2n}$. By
construction, the sets $\J_1$ and $\J_2$ are disjoint.
Let $i\in\J_2$. Then $R_i^{2n}$ is a tightly packed closed walk in the
equivalence class $\E$.
Thus, from above, there is some $j\in\J_1$ with $R_i^{2n} = R_j^{2n}$. Hence
$i=j$ and the sets $\J_1$ and $\J_2$ are
not disjoint. This gives the required contradiction. Hence $d=1$. This
completes the proof.
\end{proof}

Thus the proof of Proposition \ref{prop:coeff} is also completed.
Hence we may write $z = \alpha\sum_{i\in\I} g_i^{2n} + y$,
where $y$ is nilpotent and $\alpha \in K\setminus\{0\}$.
Moreover every tightly packed closed walk $s_2s_4\cdots s_{2n}$ in the
equivalence class $\E$ is some $R_i^{2n}$ with $i\in\I$.\medskip

So far we have only discussed non-nilpotent elements in $Z_\gr(E(\L))$
of even degree. However, the description of the non-nilpotent elements
in $Z_\gr(E(\L))$ of odd degree is easily obtained from the above as
we now explain. 

Let $z=\sum_{i=1}^l\alpha_ig_i^{2n+1}$ with $\alpha_i$ in $K\setminus
\{0\}$ for $i=1,2,\ldots,l$ be a non-nilpotent element in
$Z_\gr(E(\L))$ of odd degree $2n+1$ (with the characteristic of $K$
necessarily being equal to $2$). Observing that the degree of the
element was not used in the proofs of Lemma \ref{lemma:zero} and
Proposition \ref{prop:closedwalks}, we infer that $z=\sum_{i=1}^r z_i$
with $z_i=\sum_{g^{2n+1}_j \in\F_i}\alpha_jg^{2n+1}_j$ and each $\F_i$
is a subset of $\E_i$, and the sets $\E_i$ are distinct equivalence
classes of closed walks. Moreover, each $z_i$ is in $Z_\gr(E(\L))$,
and $z^u=\sum_{i=1}^r z_i^u$ with $z_i^u=\sum_{g^{2n+1}_j
\in\F_i}\alpha_j^u(g^{2n+1}_j)^u$. Suppose now that $z$ is
associated to the equivalence class $\E$. Then we can write
$z=\sum_{i\in \I} \alpha_ig^{2n+1}_i + y$ where $y$ is nilpotent and, 
for each $i$ in $\I$, $g^{2n+1}_i$ is non-nilpotent. Then
$z^2=\sum_{i\in \I} \alpha_i^2(g^{2n+1}_i)^2 + y^2$ is also in
$Z_\gr(E(\L))$ and is non-nilpotent of even degree! It follows from
Proposition \ref{prop:coeff} that $z=\alpha \sum_{i\in \I} g^{2n+1}_i
+ y$ where $\I=\{i: R^{2n+1}_i \text{\ is in $\E$ and $g^{2n+1}_i$ is
  not nilpotent}\}$. Moreover, it is straightforward to see, making
the obvious generalisations to the definitions, that $R^{2n+1}_i$ is
tightly packed (that is, $R^{2n+1}_i=r_{i,2}r_{i,4}\cdots
r_{i,2n}t_i$) and is both left and right $2n$-stable. In
addition, there
are paths $a_i$ and $t_i$ such that
\[R^{2n+1}_i=a_ir_{i+1,4}r_{i+1,6}\cdots r_{i+1,2n}r_{i+1,2} = 
 r_{i,2}r_{i,4}\cdots r_{i,2n}t_i\]
and
\[R^{2n+1}_i=t_{i-1}r_{i-1,2}r_{i-1,4}\cdots r_{i-1,2n}\]
with $t_ia_i=r_{i+1,2}$ (in particular $a_i=t_{i-1}$) and
$r_{i-1,2j}=r_{i+1,2j+2}$ for $j=1,2,\ldots,n-1$ and
$r_{i-1,2n}=r_{i+1,2}=r_{i,2n+1}$. 

We now summarise all these results in the following theorem.

\begin{thm}\label{thm:centre}\sloppy 
Let $z$ be a homogeneous non-nilpotent
element of degree $d$ in $Z_{\gr}(E(\L))$. Then
$$z=\sum_{i=1}^r \alpha_iz_i \mbox{ with } z_i\in Z_{\gr}(E(\L)), 
\alpha_i \in K\setminus\{0\}$$ 
where
\begin{enumerate}
\item[(i)] $\E_1, \ldots , \E_r$ are distinct equivalence classes of
closed walks,
\item[(ii)] for each $g^d$ occurring in $z_i$, the path $R^d$ is tightly
packed and lies in the equivalence class $\E_i$.
\end{enumerate}
We may write
$$z_i=\sum_{j\in \I_i}g_j^d + y_i$$
where
$\I_i = \{j \colon R_j^d \mbox{ is in } \E_i 
\mbox{ and } g_j^d \mbox{ is not nilpotent}\}$, and 
each $y_i$ is nilpotent.

If $d = 2n$, then for each $i = 1, \ldots , r$ and each $j \in \I_i$, there is some
$u \geq 1$ so that $R_j^{2n}$ stabilizes at $u$ with 
$$R_j^{2nu} = a_jr_{j+1,4}r_{j+1,6}\cdots r_{j+1, 2nu}t_j$$ and
$$R_j^{2nu} = t_{j-1}r_{j-1,2}r_{j-1,4}\cdots r_{j-1, 2nu-2}\tilde{a}_{j-1}$$
for paths $a_j$ and $\tilde{a}_j$ such that
$t_ja_j = r_{j+1, 2}$ and $\tilde{a}_{j-1}t_{j-1} = r_{j-1, 2nu}$.

If $d = 2n+1$, then for each $i,\ldots,r$ and each $j\in \I_i$ the
path $R^{2n+1}_j$ stabilizes at $1$ with 
\[R^{2n+1}_j=a_jr_{j+1,4}r_{j+1,6}\cdots r_{j+1,2n}r_{j+1,2} = 
 r_{j,2}r_{j,4}\cdots r_{j,2n}t_j\]
and
\[R^{2n+1}_j=t_{j-1}r_{j-1,2}r_{j-1,4}\cdots r_{j-1,2n}\]
with $t_ja_j=r_{j+1,2}$ (in particular $a_j=t_{j-1}$) and
$r_{j-1,2l}=r_{j+1,2l+2}$ for $l=1,2,\ldots,n-1$ and
$r_{j-1,2n}=r_{j+1,2}=r_{j,2n+1}$. 
\end{thm}

\section{Asymptotic characterization of the graded centre modulo
  nilpotents} 

We have seen in the earlier sections that non-nilpotent elements in
$Z_\gr(E(\L))$ are naturally associated to extending closed walks. For
extending closed walks the concept of being stable was introduced. In
this section we characterize non-nilpotent elements $z^u$ for any
non-nilpotent homogeneous element $z$ in $Z_\gr(E(\L))$ of even degree
and where $u$ is chosen such that the extending paths occurring in
$z^u$ are stable. We observed in section 2 that $u$ is bounded by
$\rl(\L)-1$ (Proposition \ref{prop:tpalltails}).

First we look in more detail at elements of the graded centre.  In the
next result we consider when nonzero products $g^mz \neq 0$ or $zg^m
\neq 0$ can occur, where $g^m$ is in $E(\L)$ and $z \in Z_\gr(E(\L))$.
Proposition \ref{prop:equivclass} then gives more information on the
path $R^m$ corresponding to such an element $g^m$. We keep the
notation of the previous section.

\begin{prop}\label{prop:nonzeroproducts}
With the notation of \emph{Theorem \ref{thm:centre}}, 
let $z \in Z^{2n}(E(\L))$ be non-nilpotent 
and associated to the equivalence class $\E$. Write
$z= \alpha\sum_{i\in\I}g_i^{2n} + y$ where $y = \sum_{i\in\I'}g_i^{2n}$
is nilpotent, $\alpha\in K\setminus\{0\}$.

Suppose that $g^m$ is in $E(\L)$ such that $g^mz \neq 0$ or $zg^m \neq 0$.
Then $R^m$ lies on every representative of $\E$.
\end{prop}

\begin{proof}
Suppose that $g^mz \neq 0$. Then, following Lemma \ref{lemma:zero},
there are $i, j$ with $g^mg_i^{2n} = g_j^{2n}g^m \neq 0$ and so
$R^mR_i^{2n} = R_j^{2n}R^m$.
From Lemma \ref{lem:walk}, $R^m$ is a subpath of some power of $R_i^{2n}$,
and so $R^m$ lies on $R_i^{2n}$. By the Remark following Definition
\ref{defin:liesonequiv}, we have that $R^m$ lies on every representative of $\E$.

The proof that, if $zg^m \neq 0$ then
$g^m$ lies on every representative of $\E$, is similar.
\end{proof}

We shall use tightly packed extending closed walks $R^{2n}$ to
construct elements in $Z_\gr(E(\L))$. To this end we need to know how
these paths interact with paths $R^m$ in $\R^m$ lying on $R^{2n}$. 

\begin{prop}\label{prop:equivclass}
Suppose that $R^{2n}=r_2r_4\cdots r_{2n}$ is a tightly packed, extending
closed walk, and that $R^m\in \R^{m}$ with $R^m$ lying on $R^{2n}$.
\begin{enumerate}
\item Suppose $(R^{2n})^sR^m\in \R^{2ns+m}$ with $s \geq 1$. 
  \begin{enumerate}
  \item If $m$ is even then $R^m = (R^{2n})^cr_2r_4\cdots r_{2i}$ for some 
  $c\geq 0$ and $0\leq i< n$.
  \item If\ $m$ is odd then $R^m = (R^{2n})^cr_2r_4\cdots r_{2i}p$ for some 
  $c\geq 0$, $0\leq i<n$ and path $p$ such that
  $$p\cdot t((R^{2n})^{s+c}r_2r_4\cdots r_{2i+2})=r_{2i+2}.$$
  \end{enumerate}
\item Suppose $R^m(R^{2n})^s\in \R^{2ns+m}$ with $s \geq 1$.  
  \begin{enumerate}
  \item If $m$ is even then $R^m = r_{2i}r_{2i+2}\cdots r_{2n}(R^{2n})^c$ for 
  some $c\geq 0$ and $0\leq i< n$.
  \item If $m$ is odd then $R^m = pr_{2i}r_{2i+2}\cdots r_{2n}(R^{2n})^c$ for 
  some $c\geq 0$, $0\leq i<n$ and path $p$ such that
  $$b(r_{2i-2}r_{2i}\cdots r_{2n}(R^{2n})^{s+c})\cdot p=r_{2i-2}.$$
  \end{enumerate}
\end{enumerate}
Note that if $m$ even and $i=0$ then $R^m=(R^{2n})^c$, and if $m$ odd and $i=0$ 
then $R^m = (R^{2n})^cp$, respectively $R^m = p(R^{2n})^c$.
\end{prop}

\begin{proof}  
We prove part (1) and leave part (2) to the reader.
If $m$ is even, since $R^{2n}$ is tightly packed and $R^{m}$ lies on 
$R^{2n}$, the result easily follows.
So suppose that $m$ is odd; again it is clear that
$R^m = (R^{2n})^cr_2r_4\cdots r_{2i}p$ for some path $p$. Now
$(R^{2n})^{s+c}r_2r_4\cdots r_{2i+2}$
is a prefix of $(R^{2n})^{s+c+1}$ which is in $\R^{2n(s+c+1)}$, so
by Corollary \ref{cor:tpsubpaths} we
have that $(R^{2n})^{s+c}r_2r_4\cdots r_{2i+2}\in \R^{2n(s+c)+2i+2}$.
Let $R_1 = (R^{2n})^{s+c}r_2r_4\dots r_{2i+2}$ so that 
$R_1 = (R^{2n})^sR^m\cdot t(R_1)$. Thus $p\cdot t(R_1)=r_{2i+2}$ as required.
\end{proof}

As noted before, multiplying elements in $E(\L)$ uses knowledge about
the tails and the beginnings of paths. In the next definition,
we associate not only one tail but a set of tails to an extending closed
walk.

\begin{defin} Let $R^{2n}$ be an extending closed walk. 
The {\it tail set} of
$R^{2n}$ is the set of paths $p_i$ for $0 \leq i < 2n$, such that if 
$(R^{2n})^2= R^{2n+i}q$ with $R^{2n+i}\in\R^{2n+i}$ then 
$p_i=t(R^{2n+i})$, that is, $p_i$ is the tail path of $R^{2n+i}$.  

The {\it beginning set of $R^{2n}$} is the set of paths $p_i$ for $0
\leq i < 2n$, such that if $(R^{2n})^2 = qR^{2n+i}$ with
$R^{2n+i}\in\R^{2n+i}$ then $p_i = b(R^{2n+i})$, that is, $p_i$ is the
beginning path of $R^{2n+i}$.
\end{defin}

In the next definition we look at maximal overlaps of a path with a
relation.  This is motivated by Definition \ref{defin:maximalos}, and uses 
Definition \ref{defin:overlaprelation}, which we
recall does not require the paths $P$ and $Q$ to be in $\R^2$.

\begin{defin} Let $p$ be a path in $K{\mathcal Q}$ and let $r$ be in $\R^2$.
\begin{enumerate}
\item The relation $r$ {\it maximally left overlaps} the path $p$ if 
there are paths $U, V \in K{\mathcal Q}$ with
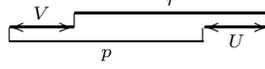
\[\xymatrix@W=0pt@M=0.3pt{%
\ar@{^{|}-^{|}}@<-1.25ex>[rrr]_p\ar@{{<}-{>}}[r]^{V} &
\ar@{_{|}-_{|}}@<1.25ex>[rrr]^{r} & & \ar@{{<}-{>}}[r]_{U} & }\] %
such that
  \begin{enumerate}
  \item[(i)] $pU = Vr$, and
  \item[(ii)] $Vr \neq PRQ$ for all non-trivial paths $Q$ and all $R \in \R^2$.
  \end{enumerate}
(Note that $V$ may be a trivial path.)
\item The path $p$ {\it maximally right overlaps} the relation $r$ if 
there are paths $U, V \in K{\mathcal Q}$ with
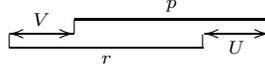
\[\xymatrix@W=0pt@M=0.3pt{%
\ar@{^{|}-^{|}}@<-1.25ex>[rrr]_r\ar@{{<}-{>}}[r]^{V} &
\ar@{_{|}-_{|}}@<1.25ex>[rrr]^{p} & & \ar@{{<}-{>}}[r]_{U} & }\] %
such that
  \begin{enumerate}
  \item[(i)] $rU = Vp$, and 
  \item[(ii)] $rU \neq PRQ$ for all non-trivial paths $P$ and all $R \in \R^2$.
  \end{enumerate}
(Note that $U$ may be a trivial path.)
\end{enumerate}
\end{defin}

The next result uses the above definitions to describe the
relations which maximally overlap a tail or a beginning path of some
$R_i^{2nu}$ where $g_i^{2n}$ occurs in a non-nilpotent element of 
$Z_\gr(E(\L))$.

\begin{prop}\label{prop:overlaprelation} 
With the notation of \emph{Theorem \ref{thm:centre}}, let $z \in 
Z^{2n}(E(\L))$ be 
non-nilpotent and associated to the equivalence class $\E$. Write $z = 
\alpha\sum_{i\in\I}g^{2n}_i + y$ with $y$ nilpotent and $\alpha \in
K\setminus\{0\}$. 

Let $r \in \R^2$ be a relation such that, for some $i \in \I$, either $r$ 
maximally left overlaps some path $p$ in the tail
set of $R_i^{2nu}$, or some path $p$ in the beginning set of
$R_i^{2nu}$ maximally right overlaps the relation $r$.
Then $r$ lies on every representative of $\E$.
\end{prop}

\begin{proof}  
Suppose the relation $r$ maximally left overlaps $p$ where $p$ is in the
tail set of $R_i^{2nu}$.  
Let $j$ be such that $0 \leq j < 2nu$, 
$(R_i^{2nu})^2 = R^{2nu+j}q$ for some path $q$ with $R^{2nu+j}\in\R^{2nu+j}$,
and $p =t(R^{2nu+j})$. Let $(r_2,r_3,\ldots,r_{2nu+j})$ be the maximal
left overlap sequence for $R^{2nu+j}$.
Since $r$ maximally left overlaps $p$ and $p=t(R^{2nu+j})$, we see that
$S=(r_2,r_3,\ldots,r_{2nu+j},r)$ is a maximal left overlap sequence.
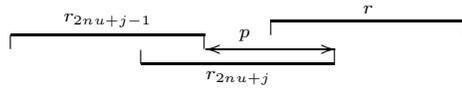
\[\xymatrix@W=0pt@M=0.3pt{%
\ar@{_{|}-_{|}}@<1.25ex>[rrr]^{r_{2nu+j-1}} & &
\ar@{^{|}-^{|}}@<-1.25ex>[rrr]_{r_{2nu+j}} & 
\ar@{<->}[rr]^(0.3){p} &
\ar@{_{|}-_{|}}@<2.5ex>[rrr]^{r} & & & & & }\]

Let $R^{2nu+j+1}$ be the element in $\R^{2nu+j+1}$ corresponding to $S$, and
let $g^{2nu+j+1}$ be the corresponding basis element of $E(\L)$.
Then $R_i^{2nu}R^{2nu+j+1}\in \R^{4nu+j+1}$ so that
$g_i^{2nu}g^{2nu+j+1}\neq 0$ and hence $z^ug^{2nu+j+1} \neq 0$. 
By Proposition \ref{prop:nonzeroproducts}, the path $R^{2nu+j+1}$ lies on
every representative of $\E$.  It follows that the last
relation in $R^{2nu+j+1}$, namely $r$, lies on every representative of $\E$. 

The other case is similar and is left to the reader.
\end{proof}

\begin{defin} A tightly packed path $r_2r_4\cdots r_{2n}$ with $r_2, r_4, 
\ldots , r_{2n} \in \R^2$ is said to be {\it relation
simple} if, whenever $i\neq j$, then $r_{2i}\neq r_{2j}$.
\end{defin}

This definition enables us to describe the minimal closed walks which
in some sense generate a non-nilpotent element in $Z_\gr(E(\L))$. 

\begin{prop}\label{prop:relationsimple}
With the notation of \emph{Theorem \ref{thm:centre}}, let $z \in 
Z^{2n}(E(\L))$ 
be non-nilpotent and associated to the equivalence class $\E$. Write $z = 
\alpha\sum_{i\in\I}g^{2n}_i + y$ with $y$ nilpotent and
$\alpha\in K\setminus\{0\}$. 

Then, for each $i \in \I$, there is a relation simple, closed walk $W_i$
such that $R_i^{2n} = W_i^{c_i}$ for some $c_i \geq 1$.  Moreover each $W_i$
lies on every representative of $\E$.
\end{prop}

\begin{proof} Let $i \in \I$ and write $R^{2n}$ for $R_i^{2n}$. 
If $R^{2n}$ is relation simple then we are done.
So suppose that $R^{2n}$ is not relation simple.
Let $R^{2n} = r_2r_4\cdots r_{2n}$.   Let $k\le n$ be the
smallest integer such that there is some $j$, with $1\le k<j\leq n$, and
$r_{2k}=r_{2j}$.  

Suppose for contradiction that $k>1$. Then 
$$r_{2j}r_{2j+2}\cdots r_{2n}r_2\cdots
r_{2j-2}\ne r_{2k}r_{2k+2}\cdots r_{2n}r_2\cdots r_{2k-2}$$
since if 
$$r_{2j}r_{2j+2}\cdots r_{2n}r_2\cdots
r_{2j-2} = r_{2k}r_{2k+2}\cdots r_{2n}r_2\cdots r_{2k-2}$$ 
then the $r_2$ on the left hand side must
equal $r_{2k-2j+2n+2}$ which contradicts $k>1$. So writing
$R^{2n}_j= r_{2j}r_{2j+2}\cdots r_{2n}r_2\cdots
r_{2j-2}$ and $R^{2n}_k=r_{2k}r_{2k+2}\cdots r_{2n}r_2\cdots
r_{2k-2}$, we have that $R^{2n}_k\neq R^{2n}_j$.
By Proposition \ref{prop:tpsequences}, we have $k, j \in \I$ and so
$g_k^{2n}$ and $g_j^{2n}$ are non-nilpotent.
Since $r_{2j} = r_{2k}$ and using the remark following Lemma
\ref{lem:relationcoeff}, this gives a contradiction to
Theorem \ref{thm:commonrelation}. Hence $k = 1$.  

Choose $j$ minimal such that $1<j\leq n$ and $r_{2j} = r_2$.
Now $R^{2n}_j$ and $R^{2n}$ both have first relation
$r_2$ so, again using Theorem \ref{thm:commonrelation}, 
we have that $R^{2n}_j = R^{2n}$. Let $W = r_2r_4\cdots r_{2j-2}$.
Then $R^{2n} = W^c$ or $R^{2n} = W^cr_2r_4\cdots r_{2s}$ for some $c\geq 1$
and $1<s<j-1$.  We show that the second case cannot occur.
For, suppose that $R^{2n} = W^cr_2r_4\cdots r_{2s}$. Then, since $R^{2n}$ is
extending, we have $r_{2s+2} = r_2$. By minimality of $j$ it follows that
$2j \leq 2s+2$ and hence $s \geq j-1$. But this contradicts the choice of
$s$. Therefore this case does not occur and $R^{2n} = W^c$.

It is now immediate from Definition \ref{defin:liesonequiv} that $W$ and 
$R^{2n}$ are in the same equivalence class, and hence $W$ lies on every
representative of the equivalence class $\E$.
Finally, by similar arguments to those given above, it follows that $W$ 
is relation simple. 
\end{proof}

A tightly packed, extending walk $R^{2n}$ need not in general be such
that the relations $r_3$, $r_5$,\ldots, $r_{2n-1}$ of the maximal left
overlap sequence form a tighly packed path $r_3r_5\cdots r_{2n-1}$,
but we show that this is asymptotically the case for the $R^{2n}$
occurring in non-nilpotent elements in $Z_\gr(E(\L))$. First we make a
precise definition describing the above phenomenon.

\begin{defin} A tightly packed, extending walk $R^{2n}$ is said to be 
{\it tightly covered} if there is some $s\geq 1$ such that 
\begin{enumerate}
\item[(i)] \sloppy if $(r_2,r_3,\ldots,r_{2ns})$ is the maximal left
overlap sequence of $(R^{2n})^s$ then
$r_{2n(s-1)-1}r_{2n(s-1)+1}\cdots r_{2ns-3}$ is a tightly packed walk,
and
\item[(ii)] if $(\tilde r_2,\tilde r_3,\ldots,\tilde r_{2ns})$ is
the maximal right overlap sequence of $(R^{2n})^{s}$ then
$\tilde{r}_3\tilde{r}_5 \cdots \tilde{r}_{2n+1}$ is a tightly packed walk.
\end{enumerate}
\end{defin}

The next result connects this notion to that of being stable.

\begin{prop}\label{prop:tightlycovered}
Suppose that $R^{2n}$ is a tightly covered closed walk with 
$t(R^{2nu})=t(R^{2n(u+1)})$ and $b(R^{2nu})=b(R^{2n(u+1)})$.
\begin{enumerate}
\item If $(r_2,r_3,\ldots , r_{2n(u+1)})$ is the maximal left overlap sequence
for $R^{2n(u+1)}$ then $r_{2nu-1}r_{2nu+1}\cdots r_{2n(u+1)-3}$ is a
tightly packed, extending walk.
\item If $(\tilde{r}_2,\tilde{r}_3,\cdots \tilde{r}_{2n(u+1)})$ is the 
maximal right overlap sequence for $R^{2n(u+1)}$ then
$\tilde{r}_3\tilde{r}_5\cdots \tilde{r}_{2n+1}$ is a
tightly packed, extending walk.
\end{enumerate}
\end{prop}

\begin{proof} 
(1) We have $t(R^{2nu})=t(R^{2n(u+1)})$ and $r_{2nu}=r_{2n(u+1)}$ so that
$r_{2nu-1}=r_{2n(u+1)-1}$.  Thus $R^{2n}$ is left $(2nu-1)$-stable. 
Hence, by  Proposition \ref{conv:p1}, for $i\geq 2nu-1$, we have 
$r_i = r_{i+2n}$.  
Since $R^{2n}$ is tightly covered, for sufficiently large $s$,
the path $r_{2ns-1}r_{2ns+1}\cdots r_{2n(s+1)-3}$ is tightly packed. 
But for $s\geq u$, $r_{2ns-1} = r_{2nu-1}$ and $r_{2ns+j} = r_{2nu+j}$ for
$0\leq j \leq 2n$. Thus $r_{2nu-1}r_{2nu+1}\cdots r_{2n(u+1)-3}$ is a tightly 
packed walk. Moreover $(r_{2nu-1},r_{2nu},\ldots,r_{2n(u+1)-3})$ is a left
overlap sequence, so by Proposition \ref{prop:tpleftoverlap}, we have that
$r_{2nu-1}r_{2nu+1}\cdots r_{2n(u+1)-3}$ is in $\R^{2n}$. Indeed
$r_{2n(u+1)-1} = r_{2nu-1}$ and so we have the following left overlap sequence
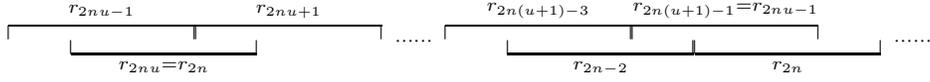
\[\xymatrix@W=0pt@M=0.3pt@C=23pt{%
\ar@{_{|}-_{|}}@<1.25ex>[rrr]^{r_{2nu-1}} & 
\ar@{^{|}-^{|}}@<-1.25ex>[rrr]_{r_{2nu}=r_{2n}} & & 
\ar@{_{|}-_{|}}@<1.25ex>[rrr]^{r_{2nu+1}} & & &
\ar@{}@<-1.25ex>[r]^{\cdots\cdots} & 
\ar@{_{|}-_{|}}@<1.25ex>[rrr]^{r_{2n(u+1)-3}} & 
\ar@{^{|}-^{|}}@<-1.25ex>[rrr]_{r_{2n-2}} & &  
\ar@{_{|}-_{|}}@<1.25ex>[rrr]^{r_{2n(u+1)-1}=r_{2nu-1}} & 
\ar@{^{|}-^{|}}@<-1.25ex>[rrr]_{r_{2n}} & & &
\ar@{}@<-1.25ex>[r]^{\cdots\cdots} & }\]
Hence from Proposition \ref{prop:tpleftoverlap}, $r_{2nu-1}r_{2nu+1}\cdots 
r_{2n(u+1)-3}$ is an extending closed walk. The result now follows.

(2) This is similar to (1).
\end{proof}

Consider all relation simple, extending closed walks in the
equivalence class $\E$.  Since the number of relations in a relation
simple walk is bounded above by the total number of (minimal)
relations generating $I$, we see that there is a finite number of
relation simple walks.  Hence the set of relation simple, extending
walks in $\E$ is finite.  We denote this set by $\Delta_{\E}$.

For each $p\in \Delta_{\E}$, we may write $p$ as a product of $s_p$
relations, that is, $p=r_2r_4\cdots r_{2s_p}$. Let $M$ be the least
common multiple of $\{s_p \mid p\in \Delta_{\E}\}$.  Then, for each
$p$, there is some $m_p \geq 1$ such that $p^{m_p}$ is the product of
$M$ relations.

By Proposition \ref{conv:p2}, for each $p\in \Delta_{\E}$, there is
$u_p \geq 1$ such that $p^{m_p}$ stabilizes at $u_{p}$.  Let $u$ be
the least common multiple of $\{u_p \mid p\in \Delta_{\E}\}$.  Then
each $p\in \Delta_{\E}$ stabilizes at $u$ and hence $p^{m_pu}$ is a
product of $Mu$ relations and has the property that
$t(p^{m_pu})=t(p^{m_p(u+1)})$ and $b(p^{m_pu}) = b(p^{m_p(u+1)})$.
Let $\Delta_{\E}^u =\{p^{m_pu} \mid p\in\Delta_{\E}\}$. For each $q\in
\Delta^u_\E$ let $g_q$ denote the corresponding element of $\G^{2N}$.

With all the previous definitions and results, we now characterize
non-nilpotent elements of sufficiently high even degree in $Z_\gr(E(\L))$.

\begin{thm}\label{thm:characterize}
With the above notation, let $\E$ be an equivalence class of closed walks in 
$\mathcal Q$ and let $W$ be a closed walk in $\E$. Let $N=Mu$. 
Then $\sum_{q\in\Delta_{\E}^u}g_q$ is in $Z_{\gr}(E(\Lambda))$ if and only if
\begin{enumerate}
\item each $q\in \Delta_{\E}^u$ is tightly covered;
\item all tightly packed, extending walks $r_2r_4\cdots r_{2N}$ that lie on 
$W$ are in $\Delta_{\E}^u$;
\item if a relation $r$ does not lie on $W$ then, for all $q\in \Delta_{\E}^u$
and all paths $p$ in the tail set of $q$, the relation $r$ does not 
maximally left overlap the path $p$;
\item if a relation $r$ does not lie on $W$ then, for all $q\in
  \Delta_{\E}^u$  
and all paths $p$ in the beginning set of $q$, the path $p$ does not maximally
right overlap the relation $r$.
\end{enumerate}
\end{thm}
\begin{proof}
Let $z =\sum_{q\in\Delta_{\E}^u}g_q$, and suppose that $z\in
Z_{\gr}(E(\Lambda))$. We show first that $q$ is tightly covered. 
Let $(r_2,r_3,\ldots,r_{4N})$ be the maximal left overlap sequence
for $q^2$. Since $q$ is left $(2N-1)$-stable, we have the overlap
sequence 
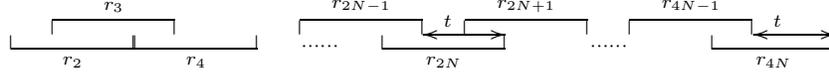
\[\xymatrix@W=0pt@M=0.3pt@C=15pt{%
\ar@{^{|}-^{|}}@<-1.25ex>[rrr]_{r_2}  &
\ar@{_{|}-_{|}}@<1.25ex>[rrr]^{r_3}   & & 
\ar@{^{|}-^{|}}@<-1.25ex>[rrr]_{r_4}  & & &
\ar@{}[rrr]_{\cdots\cdots}  &
\ar@{_{|}-_{|}}@<1.25ex>[rrr]^{r_{2N-1}} & &
\ar@{^{|}-^{|}}@<-1.25ex>[rrr]_{r_{2N}} &
\ar@{{<}-{>}}[rr]^(0.3){t} & 
\ar@{_{|}-_{|}}@<1.25ex>[rrr]^{r_{2N+1}}
& & \ar@{}[rrr]_{\cdots\cdots} & &
\ar@{_{|}-_{|}}@<1.25ex>[rrr]^{r_{4N-1}} & &
\ar@{^{|}-^{|}}@<-1.25ex>[rrr]_{r_{4N}} &
\ar@{{<}-{>}}[rr]^(0.3){t} & & }\]
\sloppy for $q^2$ with $r_{2N-1}=r_{4N-1}$ and $r_{2N}=r_{4N}$. By
Proposition \ref{relationsinmiddle}, the path $tr_{2N+2}r_{2N+4}\cdots
r_{4N-2}\tilde{a}$ where $\tilde{a}t=r_{4N}$ is the path $R^{2N}_j$
for some $j$ in $\I$ and so $g_j^{2N}$ occurs in $z$. Thus this path
is tightly packed and hence $r_{2N+1}r_{2N+3}\cdots r_{4N-1}$ is
tightly packed.  But from Proposition \ref{relationsinmiddle} again,
$t$ is the prefix of some relation and so
$\mt(r_{2N-1})=\mo(r_{2N+1})$. Hence $r_{2N-1}r_{2N+1}\cdots r_{4N-3}$
is tightly packed.  Similarly, since $q$ is right $(2N-1)$-stable, if
$(\tilde r_2,\tilde{r}_3,\ldots,\tilde r_{4N})$ is the maximal right
overlap sequence for $q^2$, then $\tilde r_{3}\tilde r_{5}\cdots\tilde
r_{2N+1}$ is tightly packed.  Thus $q$ is tightly covered.

Property (2) follows from Proposition \ref{prop:tpsequences}, and
properties (3) and (4) follow from Proposition
\ref{prop:overlaprelation}.

Now suppose that properties (1)--(4) hold, and let
$z=\sum_{q\in\Delta_{\E}^u}g_q$. Suppose first that $R^m\in \R^m$ with
$zg^m\ne 0$.  We show that $zg^m=g^mz$. First we show that $R^m$ lies
on $W$. Since $zg^m\neq 0$, we have that $qR^m$ is in $\R^{2N+m}$ for
some $q\in \Delta_E^u$. Let $(r_2,r_3,\ldots,r_{2N+m})$ be the maximal
left overlap sequence for $qR^m$.  Now $q=R^{2N-1}t(q)$ where
$R^{2N-1}\in \R^{2N-1}$ has maximal left overlap sequence
$(r_2,r_3,\ldots,r_{2N-1})$ and $r_{2N+1}$ maximally left overlaps
$t(q)$.  By hypothesis (3), since $t(q)$ is in the tail set of $q$, we
have that $r_{2N+1}$ lies on $W$.  Next, let $q^*$ be the path
associated to the maximal left overlap sequence
$(r_2,r_3,\ldots,r_{2N+1})$, so that $q^*=qt(q^*)$.  Noting that $q^*$
lies on $W$, the path $t(q^*)$ is in the tail set of $q$, and
$r_{2N+2}$ maximally left overlaps $t(q^*)$.  Thus $r_{2N+2}$ lies on
$W$.  Since $q$ is left $(2N-1)$-stable, we may continue inductively
showing that $R^m$ lies on $W$.

Suppose that $m$ is even. By Proposition \ref{prop:equivclass}, $R^m
=q^cr_2r_4\cdots r_{2i}$.  Let $q'=r_{2i+2}\cdots r_{2N}r_2\cdots
r_{2i}$, so that $q'$ is a tightly packed, extending closed walk.
Hence by (2), $g_{q'}$ occurs in $z$.  Clearly $R^mq^{\prime}=qR^m$ as
paths.  By uniqueness of the corresponding paths, $g^mz = g^mg_{q'}$ and
$zg^m=g_qg^m$. Hence $zg^m=g^mz$.

Now suppose that $m$ is odd.  Applying Proposition
\ref{prop:equivclass}, we see that $R^m =q^cr_2r_4\cdots r_{2i}p$ with
$pt(q^{c+1}r_2r_4\cdots r_{2i+2})=r_{2i+2}$.  Since $qR^m$ is in
$\R^{2N+m}$, we have that $t(qR^m)=t(R^m)=p$. Recalling that
$(r_2,r_3,\ldots, r_{4N})$ is the maximal left overlap sequence for
$q^2$, and since $q$ is tightly covered by hypothesis (1), we have
from Proposition \ref{prop:tightlycovered} that
$r_{2N-1}r_{2N+1}\cdots r_{4N-3}$ is a tightly packed, extending
walk. Now $q$ is left $(2N-1)$-stable so the path
$q'=r_{2N+2i+3}r_{2N+2i+5}\cdots r_{4N+2i+1}$ is also a tightly
packed, extending walk and
$\mt(r_{2N+2i+1})=\mt(r_{4N+2i+1})=\mo(r_{2N+2i+3})$. Thus, from
hypothesis (2), the element $g_{q'}$ occurs in $z$. 

We now show that $\mt(R^m)=\mo(q')$ so that we may consider the path
$R^mq'$. Let $s$ be the last relation in $R^m$ so $s$ is placed like
this. 
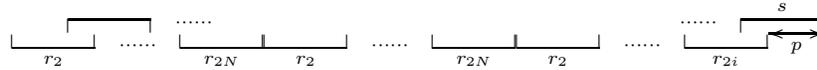
\[\xymatrix@W=0pt@M=0.3pt@C=10pt{%
\ar@{^{|}-^{|}}@<-1.25ex>[rrr]_{r_2}  & &
\ar@{_{|}-_{|}}@<1.25ex>[rrr] & 
\ar@{}[rrr]_{\cdots\cdots} & & \ar@{}@<2ex>[rrr]_{\cdots\cdots} & 
\ar@{^{|}-^{|}}@<-1.25ex>[rrr]_{r_{2N}}  & & &
\ar@{^{|}-^{|}}@<-1.25ex>[rrr]_{r_2}  & & & 
\ar@{}[rrr]_{\cdots\cdots} & & & 
\ar@{^{|}-^{|}}@<-1.25ex>[rrr]_{r_{2N}}  & & &
\ar@{^{|}-^{|}}@<-1.25ex>[rrr]_{r_2}  & & & 
\ar@{}[rrr]_{\cdots\cdots} & & \ar@{}@<2ex>[rrr]_{\cdots\cdots} & 
\ar@{^{|}-^{|}}@<-1.25ex>[rrr]_{r_{2i}}  & & 
\ar@{_{|}-_{|}}@<1.25ex>[rrr]^s & 
\ar@{{<}-{>}}[rr]_(0.5){p} & & }\]
Since $t(qR^m)=t(R^m)$, the relation $s$ is also the last relation in
$qR^m$, that is, $s=r_{2N+2i+1}$, since $q$ is $(2N-1)$-stable. Hence
$\mt(R^m)=\mt(s)=\mt(r_{2N+2i+1})=\mo(r_{2N+2i+3})=\mo(q')$. It is now
clear that $R^mq'=qR^m$ since $q'$ is tightly packed. So by uniqueness
of corresponding paths, it follows that $zg^m=g_qg^m=g^mg_{q'}=g^mz$.
Thus if $zg^m\neq 0$, then $zg^m=g^mz$.

Similarly, if $R^m\in \R^m$ with $g^mz\neq 0$ then $g^mz=zg^m$. Hence
$z$ is in $Z_\gr(E(\L))$ and this completes the proof.
\end{proof}

\section{Finite generation}

In this section we show that $Z_\gr(E(\L))/\N_Z$ is a finitely
generated commutative $K$-algebra of Krull dimension at most one. We
keep the notation of the previous sections.

We begin by showing there are only a finite number of equivalence
classes $\E$ that occur in an element of $Z_\gr(E(\L))/\N_Z$.

\begin{prop}\label{prop:finequivcl}
The number of equivalence classes that can occur in elements of the graded
centre is bounded by the number of relations in the minimal generating set
$\R^2$ that was chosen for the ideal $I$.
\end{prop}
\begin{proof}
We first consider elements of even degree. Let $z_1$ and $z_2$ be
non-nilpotent homogeneous elements in $Z_\gr(E(\L))$ of degrees $2n_1$
and $2n_2$, respectively. Suppose further that $z_i$ is associated to
the equivalence class $\E_i$. Let $R_1^{2n_1} = r_2r_4\cdots r_{2n_1}$
and $R_2^{2n_2} = \tilde{r}_2\tilde{r}_4\cdots \tilde{r}_{2n_2}$ be
paths representing non-nilpotent basis elements in $E(\L)$ occurring
in $z_1$ and $z_2$, respectively. Thus $R^{2n_1}_1$ is in $\E_1$ and
$R^{2n_2}_2$ is in $\E_2$. Suppose that $r_{2i}=\tilde r_{2j}$ for
some $i$ and $j$. By taking appropriate powers we may assume that
$n_1=n_2=n$ and that both $R_1^{2n}$ and $R_2^{2n}$ are left and right
$(2n-1)$-stable. Applying Lemma \ref{lem:relationcoeff} and Theorem
\ref{thm:commonrelation} we infer that $\E_1=\E_2$.

If $z$ is a non-nilpotent element in $Z_\gr(E(\L))$ of odd degree,
then $z^2$ is a non-nilpotent element in $Z_\gr(E(\L))$ of even
degree. Hence we can apply the above argument. It follows from this
that the number of different equivalence classes is bounded by the
number of elements in $\R^2$.
\end{proof}

The following class of examples shows that the number of equivalence
classes of closed walks associated to homogeneous non-nilpotent
elements in the graded centre of the $\Ext$ algebra of a monomial
algebra can be arbitrary large.
\begin{example}
\sloppy Choose $n$ pairs of positive integers
$\{(l_1,l_1'),(l_2,l_2'),\ldots, (l_n,l_n')\}$ where $l_i'\geq
2$. Consider the quivers $\mathcal{Q}_i=\widetilde{\mathbb{A}}_{l_i}$
with cyclic orientation.  Choose a vertex $v_1'$ of $\mathcal{Q}_1$,
two distinct vertices $v_i$ and $v_i'$ of $\mathcal{Q}_i$ for
$i=2,3,\ldots,n-1$ and finally one vertex $v_n$ of
$\mathcal{Q}_n$. Glue a copy of $\mathbb{A}_m$ with some orientation
and size $m\geq 1$ between the vertices $v_i'$ and $v_{i+1}$ for
$i=1,2,\ldots,n-1$. Call the resulting quiver $\mathcal{Q}$. Let
$\L=K\mathcal{Q}/I$ where $I$ is generated by all paths of length
$l_i'$ in subquivers $\mathcal{Q}_i$. Then $\L$ is a monomial algebra
with $n$ different equivalence classes of closed walks in the graded
centre of $E(\L)$.
\end{example}

Next we show that each equivalence class $\E$ gives rise to at most
one basis element in any given degree of $Z_\gr(E(\L))/\N_Z$.

\begin{prop}\label{prop:oneinaclass}
Let $\E$ be an equivalence class and let $n$ be a positive integer. 
\begin{enumerate}
\item[(a)] Let $V_\E$ be the subspace of $Z^n(E(\L))/\N_Z^n$ generated
  by all elements associated to $\E$. Then $\dim_K V_\E\leq 1$. 
\item[(b)] Let $r$ be the number of equivalence classes $\E$. Then
\[ \dim_K Z^n(E(\L))/\N_Z^n\leq r\]
for all $n \geq 1$.
\end{enumerate}
\end{prop}
\begin{proof}
The claim in (a) is an immediate consequence of Theorem
\ref{thm:centre}, and the statement in (b) is a direct consequence of
(a). 
\end{proof}

Now we apply the above to give further information on the
multiplication of elements associated to the same equivalence class
$\E$.

\begin{lem}\label{lem:multsameeq}
Let $z$ and $z'$ be two homogeneous non-nilpotent elements in
$Z_\gr(E(\L))$ associated to the same equivalence class $\E$. 
Then $zz'$ is a homogeneous non-nilpotent element in $Z_\gr(E(\L))$. 
\end{lem}
\begin{proof}
Let $z$ and $z'$ be two homogeneous non-nilpotent elements in
$Z_\gr(E(\L))$ associated to the same equivalence class $\E$ of
degrees $m$ and $n$, respectively. The elements $z^n$ and $(z')^m$ are
homogeneous non-nilpotent elements in $Z_\gr(E(\L))$ of degree $mn$
so, modulo nilpotence, we have $(z')^m=\alpha z^n$ for some $\alpha$ in
$K\setminus\{0\}$ by Proposition \ref{prop:oneinaclass}(a).
Furthermore,
\[(zz')^{mi}=\pm z^{mi}(z')^{mi}=\pm \alpha^i z^{(m+n)i}\neq 0\]
for all $i>0$, so that $zz'$ is a non-nilpotent element in
$Z_\gr(E(\L))$.  
\end{proof}

Let $\E$ be an equivalence class. Let $D$ be the set of positive
integers $n$ such that there is at least one non-nilpotent element $z$
in $Z_\gr(E(\L))$ of degree $n$ associated to the equivalence class
$\E$. By the above lemma the set $D$ is closed under addition and
multiplication by positive integers. Let $D = \{n_i\}_{i=1}^{\infty}$,
where $n_1 < n_2 < n_3 < \ldots.$

\sloppy To prove that there are only a finite number of generators
required for $Z_\gr(E(\L))/\N_Z$ as a $K$-algebra, we use the
following lemma. A proof can be based on the following result of James
Joseph Sylvester: Let $a$ and $b$ be positive integers. If
$\gcd(a,b)=d$, then $ab-a-b$ is the largest multiple of $d$ which
cannot be expressed in the form $ax+by$ for positive integers $x$ and
$y$.

\begin{lem}\label{lem:classical}
Let $\{a_i\}_{i=1}^{\infty}$ be a subset of
${\mathbb N}$ with $a_1 < a_2 < a_3 \ldots$. For $t\ge 1$, let
$$D(\{a_1, \ldots , a_t\}) = \{\sum_{i=1}^ta_ib_i \colon b_i \geq 0\}$$
and
$$D(\{a_i\}_{i=1}^{\infty}) = \bigcup_{t\in \mathbb N}D(\{a_1, \ldots ,
a_t\}).$$
Then there exists N such that
$$D(\{a_i\}_{i=1}^{\infty}) = D(\{a_1, a_2, \ldots , a_N\}).$$
\end{lem}

We are now ready to give the main result of this section, namely the
finite generation of $Z_\gr(E(\L))/\N_Z$ as a $K$-algebra.

\begin{thm}\label{thm:fingen}
Let $\L=K\mathcal{Q}/I$ be an indecomposable monomial algebra. Then 
$Z_\gr(E(\L))/\N_Z$ is a finitely generated commutative $K$-algebra 
of Krull dimension at most one.
\end{thm}
\begin{proof}
By Theorem \ref{thm:centre} any homogeneous non-nilpotent element $z$
of degree $n$ in $Z_\gr(E(\L))$ can be written as a sum of elements
$z_i$ in $Z_\gr(E(\L))$ of degree $n$, where each $z_i$ is associated
to a distinct equivalence class $\E_i$ for $i=1,2,\ldots,r$. In
addition the product of any two elements associated to two
different equivalence classes is zero (Proposition
\ref{prop:nonzeroproducts}). Hence it is enough to show that each
subring of $Z_\gr(E(\L))/\N_Z$ generated by elements associated to a 
fixed equivalence class is a finitely generated $K$-algebra. 

Let $\E$ be an equivalence class. Consider again the set $D$ of
positive integers $n$ such that there is at least one non-nilpotent
element $z$ in $Z_\gr(E(\L))$ of degree $n$ associated to $\E$. By
Lemma \ref{lem:classical} we can suppose that $D=D(\{
n_1,n_2,\ldots,n_t\})$ for some set $\{ n_1,n_2,\ldots,n_t\}$ of
elements in $D$. For each $n_i$ choose $z_i$ in $Z^{n_i}(E(\L))$ such
that $z_i$ is non-nilpotent. 

Suppose now that $z$ is a non-nilpotent element in $Z^n(E(\L))$
associated to $\E$. By choice of $\{ n_1,n_2,\ldots,n_t\}$ there exist
non-negative integers $\{ c_1,c_2,\ldots, c_t\}$ such that
$\sum_{i=1}^tc_in_i=n$.  Furthermore $z'=z_1^{c_1}z_2^{c_2}\cdots
z_t^{c_t}$ is a non-nilpotent element in $Z^n(E(\L))$ associated to
$\E$ by Lemma \ref{lem:multsameeq}. So from Proposition
\ref{prop:oneinaclass}, we have that $z=\alpha z'$ in
$Z_\gr(E(\L))/\N_Z$ for some $\alpha$ in $K\setminus\{0\}$. Thus
$Z_\gr(E(\L))/\N_Z$ is a finitely generated $K$-algebra.

We have already observed that $Z_\gr(E(\L))/\N_Z$ is a commutative
$K$-algebra. Moreover, by Proposition \ref{prop:oneinaclass}(b) we
have $\dim_K Z^n(E(\L))/\N_Z^n\leq r$ for some $r$ and for all $n\geq
0$. As $Z_\gr(E(\L))/\N_Z$ is a graded ring, the Krull dimension is
equal to the rate of polynomial growth of $Z^n(E(\L))/\N_Z^n$. Hence
the Krull dimension of $Z_\gr(E(\L))/\N_Z$ is at most one.
\end{proof}

We now give an example to show that $Z_\gr(E(\L))$ need not be
finitely generated, that is, the nilpotent part may be infinitely
generated.

\begin{example}
Let $\L = K{\mathcal Q}/I$ where ${\mathcal Q}$ is the quiver
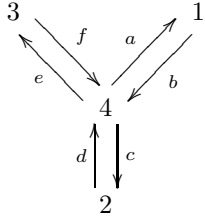
\[\xymatrix{
3\ar@<1ex>[dr]^f & & 1\ar@<1ex>[dl]^b \\
& 4\ar@<1ex>[ul]^e\ar@<1ex>[ur]^a\ar@<1ex>[d]^c & \\
& 2\ar@<1ex>[u]^d & }\]
and
$$I = \langle ab, ba, bc, cd, ef, fa \rangle.$$
By direct calculations, one shows that $Z_\gr(E(\L))$ is spanned by the
identity and by the elements $g^{2n+4}$ corresponding to the paths
$R^{2n+4}=ef(ab)^ncd$ for $n\geq 1$ as a vector space over $K$. Since
$g^{2n_1+4}g^{2n_2+4}=0$ for all pairs of positive integers $n_1$ and
$n_2$, the $K$-algebra $Z_\gr(E(\L))$ is infinitely generated. 
\end{example}

\section{The graded centre and Hochschild cohomology}

Let $\L$ be an indecomposable finite dimensional $K$-algebra.
Following \cite{SS}, let $\N$ be the ideal in $\HH^*(\L)$ generated by
the homogeneous nilpotent elements. The Hochschild cohomology ring
$\HH^*(\L)$ of $\L$ is graded commutative, so that $\HH^*(\L)/\N$ is a
commutative ring, in the same way as $Z_\gr(E(\L))/\N_Z$. As mentioned
in the introduction, it was conjectured in \cite{SS} that
$\HH^*(\L)/\N$ is always a finitely generated $K$-algebra. In this
section we show that this is indeed true whenever $\L$ is a finite
dimensional monomial algebra. Furthermore, in this case $\HH^*(\L)/\N$
has Krull dimension at most one, and $\HH^{2n}(\L)/(\N \cap
\HH^{2n}(\L))$ is isomorphic to $Z^{2n}(E(\L))/(\N_Z\cap
Z^{2n}(E(\L)))$ for $n$ sufficiently large. Note that $\N_Z\cap
Z^{2n}(E(\L))=\N_Z^{2n}$.  We end the section with a discussion on the
relationship between $\HH^*(\L)$ and $Z_\gr(E(\L))$ in general and
give an example showing that the ring homomorphism
$\overline{\varphi_{\L/\rrad}}\colon \HH^*(\L)/\N\to
Z_\gr(E(\L))/\N_Z$ is not always onto.

Let $\L=K\mathcal{Q}/I$ be an indecomposable finite dimensional
monomial algebra. From Theorem \ref{thm:fingen} we have that
$Z_\gr(E(\L))/\N_Z$ is a finitely generated $K$-algebra of Krull
dimension at most one. Moreover, using Proposition
\ref{prop:tpalltails}, there is some $u\leq \rl(\L)-1$ such that for
each generator $z+\N_Z$, the element $z$ stabilizes at $u$ and
$z^u+\N_Z$ is nonzero.  We show first that each such $z^u+\N_Z$ (or
$z^{u+1}+\N_Z$ if $u$ is odd) is in the image of the inclusion
$\overline{\varphi_{\L/\rrad}}\colon \HH^*(\L)/\N\to
Z_\gr(E(\L))/\N_Z$.

{F}rom \cite{GHZ} and Definition \ref{def:tailpaths}, recall that an
element $R_i^{2n+1}$ in $\R^{2n+1}$ may be written uniquely as
$R_k^{2n}p$ and uniquely as $qR_j^{2n}$ for some integers $j$ and $k$
and paths $p$ and $q$ in $K{\mathcal Q}$.  We consider $\L$ as a right
$\L^e$-module, and write $\otimes$ for $\otimes_K$.  {F}rom \cite{Ba},
the map $\delta^{2n+1}\colon P^{2n+1} \rightarrow P^{2n}$ in a minimal
projective $\L^e$-resolution $(P^*, \delta^*)$ of $\L$, is given by
\[\mo(R_i^{2n+1})\otimes \mt(R_i^{2n+1}) \mapsto \mo(R_k^{2n})\otimes
p - q\otimes \mt(R_j^{2n}),\]
where the first tensor lies in the
summand $\L\mo(R_k^{2n}) \otimes \mt(R_k^{2n})\L$ of $P^{2n}$ and the
second tensor lies in the summand $\L\mo(R_j^{2n}) \otimes
\mt(R_j^{2n})\L$ of $P^{2n}$.  If not specified, then it will always
be clear from the context in which summand of a projective module our
tensors lie.

\begin{prop}\label{prop:assymponto}
Keep the notation from above. Furthermore, let $\L = K{\mathcal Q}/I$
be an indecomposable monomial algebra, and let $m$ be the even integer
equal to $u$ or $u+1$. For each generator $z+\N_Z$ of
$Z_\gr(E(\L))/\N_Z$, the element $z^m$ is in the image of
$\phi_{\L/\rrad}$.
\end{prop}
\begin{proof}
Let $m$ be the even integer equal to $u$ or $u+1$. Let $z + \N_Z$ be a
non-zero homogeneous element of $Z_\gr(E(\L))/\N_Z$ associated to the
equivalence class $\E$ and suppose that $z^m$ is in degree $2n$. We
may suppose that $z^m = \sum_{i\in\I}g_i^{2n}$ where the $g_i^{2n}$
are not nilpotent, all lie in the same equivalence class, and all
$R^{2n}_i$ are left and right $(2n-1)$-stable. Moreover the $R_i^{2n}$
are tightly packed and satisfy the properties of Proposition
\ref{relationsinmiddle}.

Let $(P^*, \delta^*)$ be a minimal projective resolution of $\L$ as a
right $\L^e$-module.  Then the summands of $P^{2n}$ are indexed by the
set $\R^{2n}$.  For $R^{2n}\in\R^{2n}$, define $\chi\colon
P^{2n}\to\L$ by
\[\chi(\mo(R^{2n})\otimes \mt(R^{2n}))=
\begin{cases}
\mo(R_i^{2n}) & \text{\ if\ } R^{2n}=R^{2n}_i \text{\ with\ } i\in\I; \\
0 & \text{\ otherwise.}
\end{cases}\]
It is clear that if $\chi$ represents an element in $\HH^*(\L)$, then 
$\varphi_{\L/\rrad}(\chi)=z^m$. So, it is enough to show that $\chi$
represents an element of $\HH^{2n}(\L)$. 

Fix some $R^{2n+1} \in \R^{2n+1}$.
Consider the composite map
$$\chi \delta^{2n+1} \colon P^{2n+1} \to P^{2n} \to \L$$
and let
$\L\mo(R^{2n+1})\otimes \mt(R^{2n+1})\L$ be an arbitrary summand of
$P^{2n+1}$.  Write $R^{2n+1} = R_{i_1}^{2n}p = qR_{i_2}^{2n}$ for some
paths $p$ and $q$.

If neither $i_1$ nor $i_2$ are in $\I$ then the restriction of $\chi
\delta^{2n+1}$ to the summand $\L\mo(R^{2n+1})\otimes \mt(R^{2n+1})\L$
is zero.

Now suppose that $i_1 \in \I$, and without loss of generality suppose
that $i_1 = 1$. Then, from Proposition \ref{relationsinmiddle}, we may
write $R_1^{2n} = a_1r_{2,4}r_{2,6}\cdots r_{2,2n}t_1$ where $t_1a_1 =
r_{2,2}$.

From Proposition \ref{prop:overlaprelation} we know that $R^{2n+1}$
must lie on any walk in the equivalence class $\E$.  Now $R_1^{2n} =
R_1^{2n-1}t_1$ and there is a relation on the walk which has prefix
$t_1$, namely $r_{2,2}$. So we have that $R^{2n+1} =
R_1^{2n-1}t_1a_1 = R_1^{2n}a_1$ and hence $p=a_1$.  Then
$$\begin{array}{rcl}
R_1^{2n}a_1 & = & (a_1r_{2,4}r_{2,6}\cdots r_{2,2n}t_1)a_1\\
& = & a_1(r_{2,4}r_{2,6}\cdots r_{2,2n}r_{2,2})\\
& = & a_1R_j^{2n}
\end{array}$$
for some $j\in\I$ by Lemma \ref{lem:relationcoeff}.  Hence $q = a_1$
and the map $\delta^{2n+1}$ is given on this summand by
$$\delta^{2n+1}(\mo(R^{2n+1})\otimes \mt(R^{2n+1})) =
\mo(R_1^{2n})\otimes a_1 - a_1\otimes \mt(R_j^{2n}).$$
Hence
$$\chi\delta^{2n+1}(\mo(R^{2n+1})\otimes \mt(R^{2n+1})) =
a_1 - a_1 = 0.$$
Thus $\chi\in\HH^*(\L)$.

The final case where $i_2 \in \I$ is similar. This completes the
proof.
\end{proof}

This has the following immediate consequence.

\begin{prop}
Let $\L=K{\mathcal Q}/I$ be an indecomposable monomial algebra. Then
the ring homomorphism $\varphi_{\L/\rrad}$ induces an isomorphism
\[\HH^{2n}(\L)/(\N\cap\HH^{2n}(\L)) \cong Z^{2n}(E(\L))/(\N_Z\cap
Z^{2n}(E(\L)))\] 
for some positive integer $N$ and for all $n\geq N$. 
\end{prop}
In particular, when the characteristic of $K$ is different from two,
this shows that the map $\overline{\varphi_{\L/\rrad}}\colon
\HH^*(\L)/\N\to Z_\gr(E(\L))/\N_Z$ is an isomorphism beyond a certain
degree. If $z$ is a non-nilpotent element in $Z_\gr(E(\L))$ of odd
degree, then $z^{2i}z$ is in the $\HH^*(\L)$-submodule of
$Z_\gr(E(\L))/\N_Z$ generated by $z$ for sufficiently large $i$. Using
this it follows that $Z_\gr(E(\L))/\N_Z$ is a finitely generated
$\HH^*(\L)$-module.\medskip 

Combining the results that we have obtained so far, we prove the
main result of the paper. 
\begin{thm}
Let $\L=K\mathcal{Q}/I$ be an indecomposable monomial algebra. Then
$\HH^*(\L)/\N$ is a commutative finitely generated $K$-algebra of
Krull dimension at most one. 
\end{thm}
\begin{proof}
The map $\overline{\varphi_{\L/\rrad}}\colon \HH^*(\L)/\N\to
Z_\gr(E(\L))/\N_Z$ is an inclusion, so that 
\[\dim_K \HH^n(\L)/(\N\cap
\HH^n(\L)) \leq \dim_K Z^n(E(\L))/(\N_Z\cap Z^n(E(\L)))\leq N\] 
for some integer $N$ and for all $n\geq 0$. As $\HH^*(\L)/\N$ is a
commutative graded ring, it follows that the Krull dimension is at
most one.

If the Krull dimension is zero, there is nothing left to prove.  So
suppose that the Krull dimension of $\HH^*(\L)/\N$ is one. The Krull
dimension of $Z_\gr(E(\L))/\N_Z$ is also one and this is a finitely
generated $K$-algebra, so that by the Noether Normalisation Theorem
there exists a homogeneous element $\overline{z}=z + \N_Z$ in
$Z_\gr(E(\L))/\N_Z$ such that $Z_\gr(E(\L))/\N_Z$ is a finitely
generated module over the subring $K[\overline{z}]$. Using Proposition
\ref{prop:assymponto} there exists an element $\eta$ in $\HH^*(\L)$
such that $\overline{\varphi_{\L/\rrad}}(\eta) = \overline{z}^m$ for
some positive integer $m$. Then it easy to see that
$Z_\gr(E(\L))/\N_Z$ is a finitely generated module over the subring
$K[\overline{\eta}]$ of $\HH^*(\L)/\N$. Since $\HH^*(\L)/\N$ is a
$K[\overline{\eta}]$-submodule of $Z_\gr(E(\L))/\N_Z$, the module
$\HH^*(\L)/\N$ is a finitely generated $K[\overline{\eta}]$-module.
It is now immediate that $\HH^*(\L)/\N$ is a finitely generated
$K$-algebra.
\end{proof}
Thus the conjecture made in \cite{SS} holds for all indecomposable
finite dimensional monomial algebras $\L=K{\mathcal Q}/I$. In
\cite{GS} even stronger results are obtained for a class of algebras
that generalises Koszul monomial algebras, namely for $(D,A)$-stacked
monomial algebras. In that case the precise description and number of
generators of $\HH^*(\L)/\N$ is also found.

Now let $\L$ be any finite dimensional $K$-algebra. The proof of the
main result of this paper heavily depended on our special knowledge of
$\Im\varphi_{\L/\rrad}$ as a subring of $E(\L)$. The precise
relationship between $\Im \varphi_{\L/\rrad}$ and $E(\L)$ has been
found by Bernhard Keller \cite[Theorem 3.5]{K}.  The $\Ext$ algebra of
any module is endowed with an $\mathbb{A}_\infty$-structure (see
\cite{K2}), so using this $\mathbb{A}_\infty$-structure for $E(\L)$,
one defines a centre $Z_\infty(E(\L))$ of $E(\L)$ as an
$\mathbb{A}_\infty$-algebra.  Keller has shown that
$\Im\varphi_{\L/\rrad}=Z_\infty(E(\L))$.  Hence our results show that,
modulo nilpotence, $Z^{2n}_\gr(E(\L))$ and $Z^{2n}_\infty(E(\L))$ are
isomorphic for all $n\geq N$ for some $N$ whenever $\L$ is a finite
dimensional monomial algebra.

The $\Ext$ algebra $E(\L)$ has a trivial $\mathbb{A}_\infty$-structure
when $\L$ is a Koszul algebra, so that in this case
$Z_\infty(E(\L))=Z_\gr(E(\L))$ and $\varphi_{\L/\rrad}\colon
\HH^*(\L)\to Z_\gr(E(\L))$ is surjective. This was first independently
observed by Green-Snashall-Solberg for Koszul path algebras and for
more general Koszul algebras by Ragnar-Olaf Buchweitz \cite{Bu}, and
then subsequently generalised by Keller as referred to above \cite{K}.

The algebra $\L=K\mathcal{Q}/J^2$ is a Koszul algebra for any finite
quiver $\mathcal{Q}$ with $J$ denoting the ideal generated by all the
arrows in $\mathcal{Q}$. Here $E(\L)\cong K\mathcal{Q}$. These
algebras were considered in \cite{C}. Using the following lemma we use
the above observations to obtain some of the structural information on
the Hochschild cohomology ring of $K\mathcal{Q}/J^2$ given there.

\begin{lem}
Let $\L$ be a finite dimensional $K$-algebra such that
$\L/\rrad\otimes_K \L/\rrad$ is semisimple. Then $(\Ker
\varphi_{\L/\rrad})^{\rl(\L)} =(0)$. 
\end{lem}
\begin{proof}
Let $\eta$ be a homogeneous element in $\Ker\varphi_{\L/\rrad}$ of
degree $m$. In the proof of Proposition 4.4 in \cite{SS} it is shown
that $\Im \Omega^n_{\L^e}(\eta)\subseteq \rrad \Omega^n_{\L^e}(\L) +
\Omega^n_{\L^e}(\L) \rrad$ for any $n\geq 0$. So if $\eta'$ is a
homogeneous element in $\Ker\varphi_{\L/\rrad}$ of
degree $n$, then $\Im (\eta\Omega^m_{\L^e}(\eta'))\subseteq
\rrad^2$. Equivalently, $\Im (\eta\cdot \eta') \subseteq \rrad^2$. 
By induction we have that $\Im(\eta_1\cdot \eta_2\cdots \eta_t)\subseteq
\rrad^t$ for any set of homogeneous elements
$\{\eta_1,\eta_2,\ldots,\eta_t\}$ in
$\Ker\varphi_{\L/\rrad}$. Elements of degree zero in
$\Ker\varphi_{\L/\rrad}$ are in the radical of the centre of $\L$,
which is contained in $\rrad$. The claim follows from these
observations.  
\end{proof}

It is easy to see that $Z_\gr(K\mathcal{Q})\cong K$ whenever
$\mathcal{Q}$ is not an oriented cycle (not equal to
$\widetilde{\mathbb{A}}_n$ with cyclic orientation for any $n$). Again
let $\L=K\mathcal{Q}/J^2$.  Hence if $\mathcal{Q}$ is not an oriented
cycle, then $\Ker \varphi_{\L/\rrad} = \HH^{\geq 1}(\L)\cup \rad
\HH^0(\L)$. Using the above lemma we obtain that any product of two
elements in $\HH^{\geq 1}(\L)\cup \rad \HH^0(\L)$ is zero. For further
results see \cite{C}.

We end this paper with two examples. The first example shows that the
map $\overline{\varphi}_{\L/\rrad}\colon \HH^*(\L)/\N\to
Z_\gr(E(\L))/\N_Z$ can be an isomorphism even though the stability
condition in Proposition \ref{prop:assymponto} is not satisfied. The
last example provides a case where the map
$\overline{\varphi}_{\L/\rrad}$ is not onto, or equivalently, is not
an isomorphism.

\begin{example}
Consider the Example \ref{ex:gradedcentre} again. The element
$z=g^4_1+g^4_2+g^4_3+g^4_4$ is non-nilpotent in $Z_\gr(E(\L))$ and
$Z_\gr(E(\L))/\N_Z\cong K[\overline{z}]$. Moreover $z$ has degree
$4$, and it is left and right $7$-stable. So the point $u$ where $z$
stabilizes is $2$, and therefore $\overline{z}^2$ is in the image of
$\overline{\varphi}_{\L/\rrad}$ by Proposition
\ref{prop:assymponto}. However we claim that $\overline{z}$ is in the
image of $\overline{\varphi}_{\L/\rrad}$. Let $P^*\colon \cdots\to
P^n\to P^{n-1}\to \cdots \to P^1\to P^0\to \L\to 0$ be a minimal
projective $\L^e$-projective resolution of $\L$. Define $\chi\colon
P^4\to \L$ by 
\[\chi(\mo(R^4_j)\otimes \mt(R^4_j))=
\begin{cases}
\mo(R^4_j)  & \text{\ for\ } j=1,2,3,4;\\
a_2 & \text{\ for\ } j=5.
\end{cases}\]
Then $\chi$ is in $\HH^4(\L)$ and $\varphi_{\L/\rrad}(\chi)=z$, so
that $\overline{\varphi_{\L/\rrad}}\colon \HH^*(\L)/\N\to
Z_\gr(E(\L))/\N_Z$ is an isomorphism. 
\end{example}

\begin{example}
Let $\L = K{\mathcal Q}/I$ where ${\mathcal Q}$ is the quiver
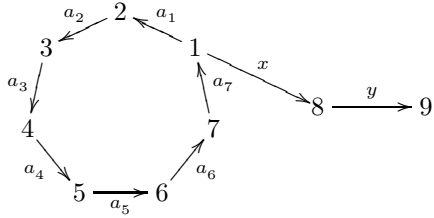
\[
\def\alphanum{\ifcase\xypolynode\or 1\or 2
\or 3\or 4\or 5\or 6\or 7\fi}
\xy/r3pc/:
\xymatrix@C=30pt{{{\xypolygon7"A"{~<<{@{}}~><{}
~*{\alphanum}
~>>{_{a_\xypolynode^{}}}}}} & & 8\ar@{}[l]^(0){}="B"\ar[r]^y &
9\ar@{{}-{>}}^x "A1";"B"}
\endxy\]
and
\[I = \langle a_1a_2a_3a_4, a_2a_3a_4a_5, a_4a_5a_6, a_5a_6a_7,
a_6a_7a_1a_2, a_6a_7xy, a_7a_1a_2a_3\rangle.\]
\sloppy It is easy to see that 
$R^4_1 = a_1a_2a_3a_4a_5a_6a_7$, 
$R^4_2 = a_4a_5a_6a_7a_1a_2a_3$,
$R^4_3 = a_5a_6a_7a_1a_2a_3a_4$, and 
$R^4_4 = a_7a_1a_2a_3a_4a_5a_6$
are all in $\R^4$. Moreover $z=g^4_1+g^4_2+g^4_3+g^4_4$ is in
$Z_\gr(E(\L))$ in degree $4$. Since one can check that
$\HH^4(\L)=(0)$, the induced map $\overline{\varphi_{\L/\rrad}}\colon
\HH^*(\L)/\N\to Z_\gr(E(\L))/\N_Z$ is not onto. To this end note that
$g^4_1$ is not stable.
\end{example}


\begin{thebibliography}{99}
\bibitem{Ba} Bardzell, M. J., {\em The alternating syzygy behavior of
    monomial algebras}, J.\ Algebra {\bf 188} (1997), 69-89.

\bibitem{Bu} Buchweitz, R.-O., \emph{Hochschild cohomology of Koszul
    algebras}, Talk at the Conference on Representation Theory,
  Canberra, July 2003.

\bibitem{C} Cibils, C., \emph{Hochschild cohomology algebra of radical
    square zero algebras}, Algebras and modules, II (Geiranger, 1996),
  93--101, CMS Conf.\ Proc., 24, Amer.\ Math.\ Soc., Providence, RI,
  1998.

\bibitem{E} Evens, L., {\em The cohomology ring of a finite group},
  Trans.\ Amer.\ Math.\ Soc.\ \textbf{101} (1961), 224--239.

\bibitem{FS} Friedlander, E. M.\ and Suslin, A., \emph{Cohomology of
    finite group schemes over a field}, Invent.\ Math.\ \textbf{127}
  (1997), no.\ 2, 209--270.

\bibitem{GHZ} Green, E.\ L., Happel, D. and Zacharia, D.,
  \emph{Projective resolutions over Artin algebras with zero
    relations}, Illinois J. Math.\ \textbf{29} (1985), no.\ 1,
  180--190.

\bibitem{Go} Golod, E., \emph{The cohomology ring of a finite
$p$-group}, Dokl.\ Akad.\ Nauk SSSR, 125 (1959) 703--706.

\bibitem{GS} Green, E.\ L.\ and Snashall, N., \emph{The Hochschild
cohomology ring modulo nilpotence of a stacked monomial algebra},
Preprint 2004. 

\bibitem{GSS} Green, E.\ L., Snashall, N.\ and Solberg, \O., \emph{The
    Hochschild cohomology ring of a selfinjective algebra of finite
    representation type}, Proc.\ Amer.\ Math.\ Soc.\ \textbf{131}
  (2003), no.\ 11, 3387--3393.
  
\bibitem{GZ} Green, E.\ L.\ and Zacharia, D., \emph{The cohomology ring
    of a monomial algebra}, Manuscripta Math.\ \textbf{85} (1994),
  no.\ 1, 11--23.
  
\bibitem{K2} Keller, B., \emph{$A$-infinity algebras in representation
    theory}, Proceedings of ICRA IX, Beijing 2000.
 
\bibitem{K} \bysame, \emph{Derived invariance of higher structures on
the Hochschild complex}, Preprint 2003.

\bibitem{SS} Snashall, N.\ and Solberg, \O., {\em Support varieties and
Hochschild cohomology rings}, to appear in Proc. London Math.\ Soc.

\bibitem{V} Venkov, B.\ B., \emph{Cohomology algebras for some
classifying spaces}, Dokl.\ Akad.\ Nauk SSSR \textbf{127} (1959), 
943--944.
\end{thebibliography}
\end{document}